\documentclass{amsart}
\usepackage{amsfonts}
\usepackage{amssymb}
\usepackage{amsmath}
\usepackage{graphicx}
\usepackage{float}
\usepackage{subfig}
\usepackage{svg} 
\usepackage{tikz}
\usepackage[spanish]{babel}
\usepackage{hyperref}
\usepackage[spanish]{cleveref}
\usepackage{enumitem,array,tabularx}
\usepackage{verbatim}
\usepackage{listings}
\usepackage{matlab-prettifier} 

\usepackage{csquotes}

\newtheorem{theorem}{Teorema}[section]
\newtheorem{corollary}[theorem]{Corolario}
\newtheorem{definition}[theorem]{Definici\'on}
\newtheorem{proposition}[theorem]{Proposici\'on}
\newtheorem{remark}[theorem]{Remark}
\newtheorem{example}[theorem]{Ejemplo}
\newtheorem{lemma}[theorem]{Lema}

\newcommand{\teorema}[2]{\begin{theorem}[#1] #2\end{theorem}}

\newcommand{\dfe}[2]{\begin{definition}[#1]#2\end{definition}}
\newcommand{\proposicion}[2]{\begin{proposition}[#1]#2\end{proposition}}
\newcommand{\lema}[2]{\begin{lemma}[#1]#2\end{lemma}}
\newcommand{\nota}[2]{\begin{remark}[#1]#2\end{remark}}

\newcommand{\demostracion}[1]{\begin{proof}#1\end{proof}}

\crefname{section}{secci\'on}{secciones}
\Crefname{section}{Secci\'on}{Secciones}

\newcommand{\Z}{\mathbb{Z}}

\newcommand{\R}{\mathbb{R}}

\begin{document}

\title{Nudos y Superficies}
%
\author[Y. Arciniegas]{Yoseth Arciniegas Barreto}
\email{ysarciniegasb@udistrital.edu.co}
\author[N. Bermudez]{Nicol Bermudez Bohorquez}
\email{Nicolbermudez777@gmail.com}
\author[J. Pinz\'on-Caicedo]{Juanita Pinz\'on-Caicedo}
\thanks{The third author was partially supported by Simons Collaboration grant 712377.}
\address {Department of Mathematics, University of Notre Dame, Notre Dame, IN 46556}
\email{jpinzonc@nd.edu}
\author[L. Rozo]{Luisa Rozo Posada}	
\email{lfrozopo@unal.edu.co}

%
%
%
%
%

\begin{abstract} Estas notas son una introducci\'on a la teor\'ia de nudos desde la perspectiva de las superficies. Las notas abordan conceptos fundamentales tales como isotop\'ias, movimientos de Reidemeister, nudos t\'oricos y superficies orientables, conexas y con una componente conexa. De igual forma se exponen invariantes de nudos definidos a trav\'es de las superficies de Seifert y sus matrices asociadas, incluyendo el 3-g\'enero, el polinomio de Alexander y la signatura. Finalmente, a partir de la noci\'on de concordancia y la operaci\'on de suma conexa, se introduce una estructura de grupo en el conjunto de nudos.
		
\end{abstract}

\maketitle
%
%

\section{Introducci\'on}\label{intro}
La \textbf{topolog\'ia} es una rama de las matem\'aticas que estudia aquellas propiedades de los objetos geom\'etricos que se mantienen bajo alguna relaci\'on de equivalencia, tal es el caso de nuestros objetos de estudio, los \textbf{nudos}. Estos, vistos desde la topolog\'ia, resultan ser objetos interesantes que nos permiten movilizarnos entre diversos espacios y dimensiones, as\'i mismo preguntarnos por sus deformaciones que permitan equivalencias entre ellos y por ende definir invariantes que se mantengan bajo estas clasificaciones.

La teor\'ia de nudos naci\'o como un intento de comprender el universo. Alrededor de 1867, cuando la comunidad cient\'ifica quer\'ia poder explicar los diferentes tipos de materia, el matem\'atico y f\'isico escoc\'es Peter Guthrie Tait le mostr\'o a su amigo Lord Kelvin (Sir William Thomson) un artefacto que al utilizarlo generaba anillos de humo. Thomson quedo fascinado con tal demostraci\'on por la belleza de los anillos y se pregunt\'o si como los anillos eran v\'ortices en el aire, los \'atomos eran anillos de v\'ortice enredados en el \'eter, un medio en el cual cre\'ian los f\'isicos que yac\'ia la energ\'ia.

La teor\'ia de v\'ortices gano popularidad e inspiro a Tait a comenzar a tabular todo los nudos, para as\'i generar algo equivalente a una tabla de elementos. Sin embargo a finales de la d\'ecada de 1880, Thomson fue abandonando su teor\'ia del v\'ortice, pero para ese entonces Tait ya estaba cautivado por la belleza matem\'atica que ocultaban sus nudos, as\'i que continuo con su proyecto de tabulaci\'on. En el proceso estableci\'o el campo matem\'atico de la teor\'ia de los nudos. \cite{knot_history} 

\'Estas notas son el resultado de el mini-curso ``Nudos y Superficies'' parte del primer ECOGyT que se llev\'o a cabo en Bogot\'a en Julio de 2024. Los resultados que aqu\'i presentamos no son originales, y muchos de los teoremas se enuncian pero no se demuestran. No obstante, intentamos incluir referencias que contienen demostraciones completas y claras.

\section{Superficies}\label{superficies}
Para poder adentrarnos al mundo de los nudos ser\'a necesario en primer lugar familiarizarnos con algunos espacios topol\'ogicos fundamentales. Para empezar, consideramos espacios conexos de dimensi\'on $0$ y notamos que \'estos son conjuntos de un solo punto. En dimensi\'on $1$ tenemos un conjunto muy especial: el c\'irculo de radio uno, definido como 
 \[S^1=\{(x,y)\in \R^2 \mid x^2+y^2=1\}. \]
Dicho conjunto se caracteriza por ser compacto, conexo y por no tener frontera. Vemos entonces que en cada una de \'estas dos dimensiones, tenemos un \'unico espacio cerrado y conexo. Tal no es el caso en las siguientes dimensiones. Los espacios de dimensi\'on 2 son conocidos como superficies, y el primer ejemplo es  la 2-esfera en $\mathbb{R}^3$:
\[S^{2}=\{(x,y,z)\in\mathbb{R}^3 \mid x^2+y^2+z^3=1\}.\]
Como el punto y el c\'irculo, $S^2$ es compacta y sin frontera (i.e. cerrado), y consiste de una componente conexa. Sin embargo, como ver\'emos a continuaci\'on, no es el \'unico objeto de dimensi\'on 2 con \'estos atributos. 

Dentro de las superficies estaremos interesados en algunas que cumplan con dos caracter\'isticas especiales: ser orientable (\cref{orientable}) y cerrada (\cref{cerrada}). 

\dfe{Superficie cerrada}{\label{cerrada}
     Un espacio topol\'ogico $S$ es una superficie cerrada si para todo punto $p\in S$, existe un abierto $U\subseteq S$ tal que $p\in U$  y $U$ es homeomorfo a una $\epsilon$-vecindad euclideana en $\R^2$. 
}

N\'otese que una superficie cerrada es en particular compacta y no tiene frontera. Por ejemplo, el disco $D^2=\{(x,y)\in \R^2 \mid x^2+y^2\leq 1\}$ (con la topolog\'ia del subespacio) es compacto, pero tiene frontera por lo que no es una superficie cerrada \footnote{Es importante no confundir la noci\'on de superficie cerrada \Cref{cerrada} con la noci\'on de subconjunto cerrado. El disco es un subconjunto cerrado de $\R^2$, pero no es una superficie cerrada.}. Para verlo, note que si $p\in D^2$ tiene norma exactamente $1$, la intersecci\'on de un disco de radio $\epsilon>0$ centrado en $p$ y el disco $D^2$ no es homeomorfo a un abierto de $\R^2$, pues contiene todos los puntos de un arco del c\'irculo $S^1$.\\

Para entender la noci\'on de orientabilidad, introducimos primero el ejemplo m\'as concreto de un espacio topol\'ogico no-orientable: la cinta de M\"obius $M$, la cual se puede  construir a partir de una tira rectangular de papel a la que se da media vuelta en el sentido `largo', para luego pegar los dos extremos cortos. \'Esta construcci\'on se puede considerar de manera abstracta como el espacio cociente del rect\'angulo $[0,10]\times [0,1]$ con los bordes verticales $\{0,10\}\times [0,1]$ identificados por la relaci\'on de equivalencia $(0,y)\sim (10,1-y)$ para $0\leq y\leq 1$.\\

\begin{figure}[h]
    \centering
    \includegraphics[width=0.675\linewidth]{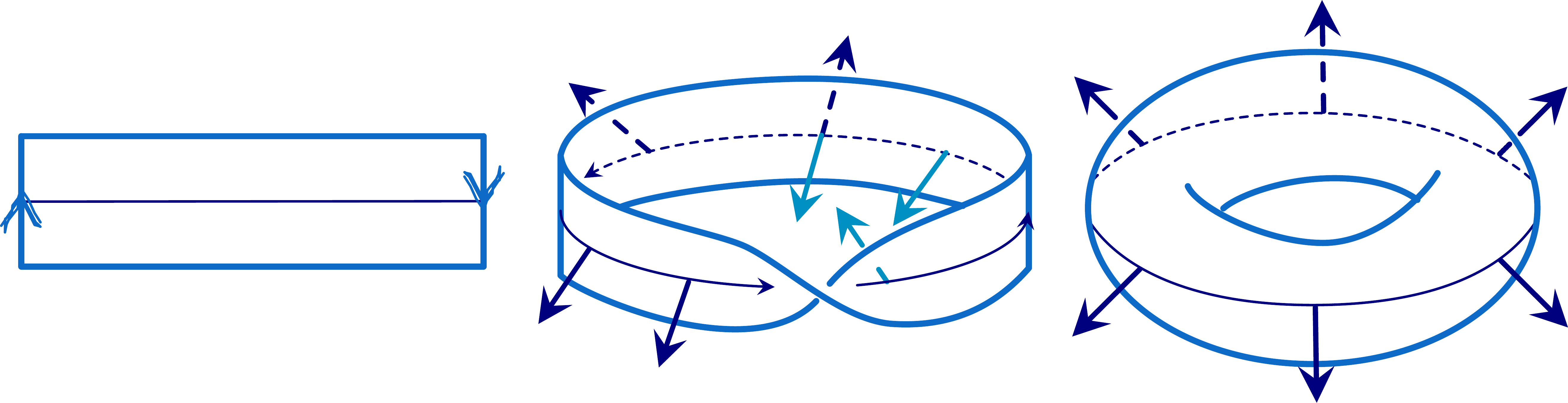}
    \caption{Una superficie no orientable y una orientable}\label{no_si_orientable}
\end{figure}

Podemos ver que $M$ no es una superficie notable: si se parte con una pareja de ejes perpendiculares orientados, al desplazarse paralelamente a lo largo de la cinta, se volverá al punto de partida con la orientación invertida. La \Cref{no_si_orientable} ilustra como la \textit{`torcerdura'} de $M$ genera un problema en la orientaci\'on. En cambio el toro si resulta ser una superficie orientable, ya que como se observa es posible encontrar ejes perpendiculares que var\'ian continuamente sobre el toro, y cuya direcci\'on no cambia.

\dfe{Superficie orientable}{\label{orientable}
Una superficie $S$ se dice no-orientable si existe una banda de M\"obius completamente contenida en $S$.}

\subsection{Presentaciones poligonales}
As\'i como la banda de M\"obius se puede describir como un espacio cociente de un rect\'angulo, toda superficie puede ser construida como el espacio cociente de un conjunto de pol\'igonos. \'Esta construcci\'on se conoce como \textit{`presentación poligonal'} y consiste de 
\begin{enumerate}
\item un conjunto de polígonos planos,
\item una orientaci\'on de las aristas de todos los pol\'igonos, y
\item una relaci\'on de equivalencia que empareja las aristas.
\end{enumerate}

\begin{figure}[h]
    \centering
    \includegraphics[width=0.5\linewidth]{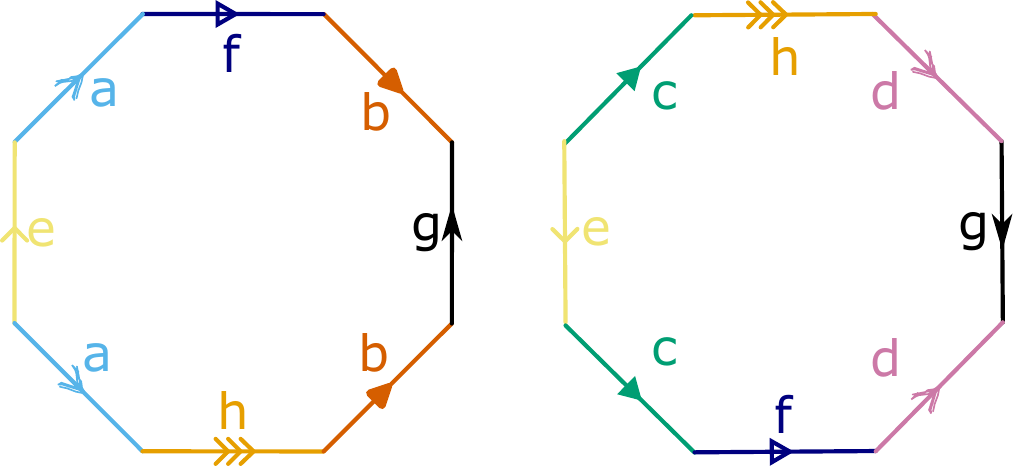}
    \caption{Una representaci\'on poligonal.}
    \label{pol-2}
\end{figure}

N\'otese que si la relaci\'on de equivalencia incluye todas las aristas,  la presentaci\'on poligonal representa una superficie cerrada. Asimismo, si el emparejamiento de las aristas conecta todos los polígonos, entonces la presentaci\'on poligonal representa una superficie conexa (de hecho en \'este caso se puede considerar una presentaci\'on poligonal con un \'unico pol\'igono). \Cref{pol-2} muestra una presentación poligonal con dos pol\'igonos.

\begin{figure}[h]
    \centering
    \def\svgwidth{0.3\textwidth}
    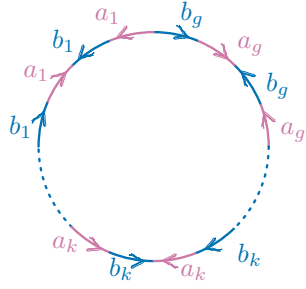
    \caption{Representac\'ion poligonal can\'onica de la suma conexa de $g$ toros. }
    \label{pol-can}
\end{figure}

\begin{figure}[h!]
    \centering
    \includegraphics[width=0.4\linewidth]{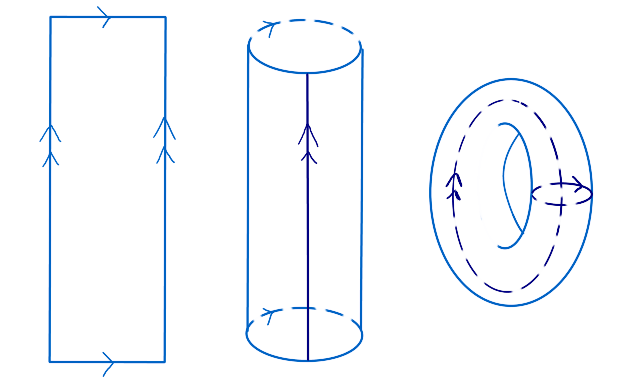}
    \caption{Pol\'igono fundamental del toro}
    \label{identificacion_toro}
\end{figure}

Cabe anotar que las representaciones poligonales no son \'unicas. De hecho, una superficie $S$ admite un n\'umero infinito de presentaciones poligonales. Con miras en la clasificaci\'on de superficies cerradas y orientables, para cada entero positivo $g\geq 1$ se define una\textit{`presentación poligonal can\'onica'} que consiste de un pol\'igono regular de $4g$ lados con el emparejamiento de sus aristas como lo muestra \Cref{pol-can}. Por ejemplo si tomamos el cuadrado ($g=1$) con la orientaci\'on e identificaci\'on que se muestra en la \Cref{identificacion_toro}, obtenemos la superficie del toro. 

\begin{figure}[h]
    \centering
    \text{Modificaci\'on 1: Eliminaci\'on de aristas contiguas y emparejadas, y de el v\'ertice que las une.}\\[1em]
    \includegraphics[height=0.15\linewidth]{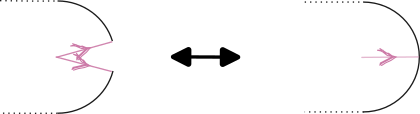}\\[1em]
    \text{Modificaci\'on 2: Divisi\'on de una cara y creaci\'on de dos aristas emparejadas.}\\[1em]
    \includegraphics[height=0.15\linewidth]{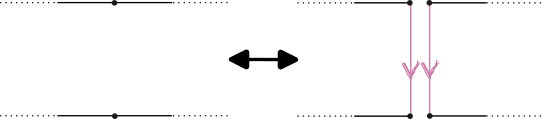}\\[1em]
    \text{Modificaci\'on 3: Divisi\'on de un par de aristas identificadas y creaci\'on de un v\'ertice.}\\[1em]
    \includegraphics[height=0.15\linewidth]{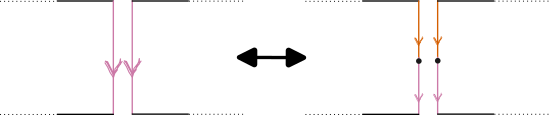}\\[1em]
    \caption{Modificaciones de presentaciones poligonales.}
    \label{modificaciones}
\end{figure}

Se dice que dos presentaciones $P$ y $P'$ son equivalentes si existe una secuencia de modificaciones \Cref{modificaciones} que transforma una en la otra. Puesto que ninguna de estas modificaciones cambia realmente la superficie subyacente, se deduce que las superficies representadas por presentaciones poligonales equivalentes son homeomorfas. De manera an\'aloga, si dos superficies $S$ y $S'$ son homeomorfas, entonces el homeomorfismo transforma una presentación poligonal de una en una presentación poligonal de la otra. Por ejemplo, \Cref{pol-2} y \Cref{pol-can} con $g=2$ son equivalentes.\\

\subsection{Caracter\'istica de Euler}
La construcci\'on de superficies como espacios cocientes de pol\'igonos nos permite entonces crear superficies a partir de v\'ertices, aristas y caras. Con esta definici\'on se logra definir un invariante de las superficies conocido como \textit{Caracter\'istica de Euler}.

\dfe{Caracter\'istica de Euler}{
Sea $M$ una superficie. La \textbf{caracter\'istica de Euler} asociada a $M$ es:
\begin{equation}\label{euler}
\chi(M)=V-A+C
\end{equation}
donde $V$, $A$ y $C$ representan n\'umero de v\'ertices, aristas y caras respectivamente.
}

\begin{figure}[b]
    \centering
    \includegraphics[width=0.3\linewidth]{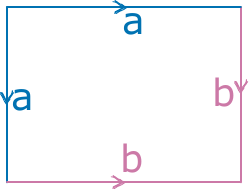}
    \caption{Esfera $S^2$ con v\'ertices, aristas y caras}
    \label{esfera}
\end{figure}

Por ejemplo, la esfera $S^2$ se puede formar a partir de 3 v\'ertices, 2 aristas y 1 cara (ver \Cref{esfera}), luego su caracter\'istica de Euler es $\chi(S^2)=3-2+1=2$. Asimismo, \Cref{identificacion_toro} muestra una descomposici\'on del toro $\mathbb{T}^2$ a partir de 1 v\'ertice, 2 aristas y 1 cara, entonces $\chi(\mathbb{T}^2)=1-2+1=0$. En general, la superficie $\Sigma_g$ constru\'ida como la suma conexa de $g$ toros admite una descomposisi\'on con 1 v\'ertice, $2g$ aristas y 1 cara, por lo que 
\begin{equation}\label{euler-g}
\chi(\Sigma_g)=1-2g+1=2-2g.
\end{equation}

En t\'erminos de las modificaciones de presentaciones poligonales que se muestran en \Cref{modificaciones}, n\'otese que las modificaciones preservan la caracter\'istica de Euler.  El \Cref{euler-modificaciones} muestra los cambios precisos en el n\'umero de v\'ertices, aristas y caras. Dejamos los c\'alculos precisos para el lector interesado.

\begin{table}[H]
\renewcommand{\arraystretch}{1.5}
\begin{tabular}{| l |c|c|c|c|c|c|}\hline
&\multicolumn{2}{|c|}{Modificaci\'on 1}&\multicolumn{2}{|c|}{Modificaci\'on 2}&\multicolumn{2}{|c|}{Modificaci\'on 3}\\\hline
& $P_\text{izq}$ & $P_\text{der}$& $P_\text{izq}$ & $P_\text{der}$& $P_\text{izq}$ & $P_\text{der}$\\\hline
V\'ertices	& $V$ & $V-1$		& $V$& $V$		& $V$& $V+1$\\
Aristas 	& $A$ & $A-1$ 		& $A$ & $A+1$ 	& $A$ & $A+1$ \\
Caras 	& $C$& $C$		& $C$& $C+1$		& $C$& $C$\\\hline
\end{tabular}
\caption{Cambios en el n\'umero de v\'ertices, aristas y caras luego de realizar las modificaciones de \Cref{modificaciones}. Se denota por $P_\text{izq}$ y $P_\text{der}$ las presentaciones poligonales que aparecen en la parte izquierda y derecha de la imagen que ilustra cada modificaci\'on.}
\label{euler-modificaciones}
\end{table}

\subsection{Clasificaci\'on de superficies cerradas y orientables}
La existencia de presentaciones poligonales para las superficies cerradas y la caracter\'istica de Euler son los elementos fundamentales en la demostraci\'on del siguiente teorema: 

\teorema{Homeomorfismo de superficies}{Una superficie orientable, cerrada  y conexa es homeomorfa la esfera o a la suma conexa de $g$ toros con $g\geq 1$. \label{homeomorfismo}}

El n\'umero $g$ de toros de la suma conexa se conoce como el \textbf{genero de la superficie} $S$ y se representa como $g(S)$. En t\'erminos de la demostraci\'on, \Cref{homeomorfismo} se obtiene como consecuencia del siguiente resultado: 

\teorema{Clasificaci\'on de superficies cerradas}{\label{clasificacion} 
Sean $S_1$ y $S_2$ superficies cerradas y orientables. Las siguientes afirmaciones son equivalentes:
\begin{enumerate}[label=(\roman*)]
\item Existe un homeomorfismo $h:S_1\to S_2$.
\item $S_1$ y $S_2$ admiten presentaciones poligonales equivalentes.
\item $\chi(S_1)=\chi(S_2)$
\end{enumerate}
}

Los detalles de la demostraci\'on se encuentran en \cite[Section 3.2]{stahl-stenson}. En t\'erminos generales, el resultado se obtiene luego de establecer los siguientes lemmas:

\begin{enumerate}[label=\Alph*.]
\item Dos superficies son homeomorfas si y s\'olo si tienen presentaciones poligonales equivalentes.
\item Cualquier presentaci\'on poligonal es equivalente a una de las presentaciones poligonales can\'onicas.
\item Si dos presentaciones poligonales son equivalentes, entonces tienen la misma caracter\'istica de Euler.
\end{enumerate}

\section{Nudos}\label{nudos}
\subsection{Generalidades}

A lo largo de esta secci\'on se desarrollar\'an algunos conceptos claves para la compresi\'on de la teor\'ia de nudos como encajamientos, isotop\'ia, equivalencia, etc...\\

Como un primer acercamiento visual nos podemos imaginar un nudo como tomar una cuerda \textit{\enquote*{enredarlarla}} y unir sus extremos. Visto desde un punto de vista matem\'atico  un nudo $K$ se puede entender como incrustar el lazo cerrado $S^1$ en el espacio 3-dimensional. Veamos esta definici\'on un poco m\'as formal.

\dfe{Nudo}{
    Un nudo $K$ en $\mathbb{R}^3$ es una curva cerrada y suave. Es decir existe una funci\'on 
    $$K:S^1\rightarrow\mathbb{R}^3$$
    suave y biyectiva en su imagen (encajamiento).
    }
    
\begin{figure}[H]
        \centering
        \subfloat[\centering ]{{\includegraphics[height=0.1\textheight]{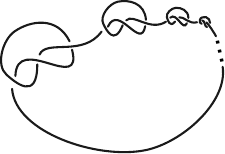} }}%
        \qquad
        \subfloat[\centering ]{{\includegraphics[height=0.1\textheight]{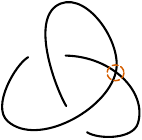} }}%
        \caption{Estos casos se evitan con suavidad y sobreyectividad}%
        \label{suave y sobre}%
    \end{figure}

La condici\'on de suavidad garantiza que no exista una sucesi\'on de nudos m\'as peque\~nos (ver Figura \ref{suave y sobre}a), por otro lado la sobreyectividad sobre la imagen es necesaria para que en la imagen no hayan cortes, es decir, se puede decidir siempre cu\'al parte del lazo va por arriba y cu\'al parte por debajo  (ver Figura \ref{suave y sobre}b).\\
    

\begin{figure}[H]
        \centering
        \includegraphics[width=0.5\linewidth]{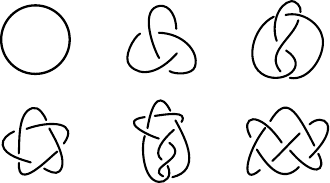}
        \caption{Ejemplos de nudos}
        \label{nudos-ejemplos}
    \end{figure}

    Este encajamiento con el que se define un nudo se puede representar gr\'aficamente apartir de un diagrama, el cual se puede interpretar como una proyecci\'on del nudo en el plano. A continuaci\'on se muestran algunos ejemplos de nudos, m\'as espec\'ificamente diagramas de diferentes nudos. Uno de esos nudos es conocido como el nudo trivial, este nudo tiene como caracter\'istica que no est\'a anudado
    (ver Figura \ref{nudos-ejemplos}).

\begin{figure}[H]
        \centering
        \includegraphics[width=0.4\linewidth]{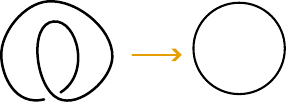}
        \caption{Deformar un nudo en el trivial}
        \label{desanudar}
    \end{figure}

    Ahora observemos el diagrama  (ver Figura \ref{desanudar}), es f\'acil observar que podemos deformar este nudo de forma suave y obtener el nudo trivial, es natural entonces preguntarnos si dado un nudo, este se puede deformar en otro nudo que en esencia sean el mismo. En otras palabras, queremos determinar cu\'ando dos nudos son equivalentes.
    
\dfe{Equivalencia de nudos}{
    Dos nudos $K_0, K_1$ nudos, son equivalentes si existe una funci\'on $F:[0,1]\times\mathbb{R}^3\longrightarrow\mathbb{R}^3$ tal que:
    \begin{enumerate}
        \item $f_t(x)=F(t,x)$ es un difeomorfismo para todo $t\in[0,1]$.
        \item $f_0(x)=x$ for all $x\in \R^3$.
        \item $f_1(K_0)=(K_1)$.
    \end{enumerate}
    }    


Esta definci\'on nos indica que comenzamos con los nudos originales en el espacio y vamos deformando $K_0$ de forma suave y que cada paso se pueda devolver de forma suave, as\'i hasta que en el \'ultimo movimiento se obtenga el nudo $K_1$.\\

\begin{figure}[h!]
            \centering
            \includegraphics[width=0.6\linewidth]{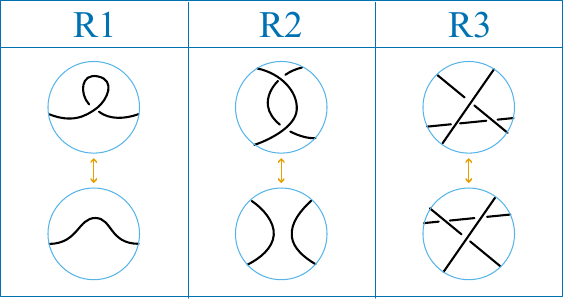}
            \caption{Movimientos de Reidemeister}
            \label{tabla reidemeister}
        \end{figure}

Sin embargo, Redeimeister en 1927 propone tres movimientos \textit{locales} que permiten transformar nudos de una forma adecuada como se solicitaba. Estos movimientos se muestran en la  tabla expuesta en la Figura \ref{tabla reidemeister}.\\

\begin{figure}[H]
            \centering
            \includegraphics[width=0.7\linewidth]{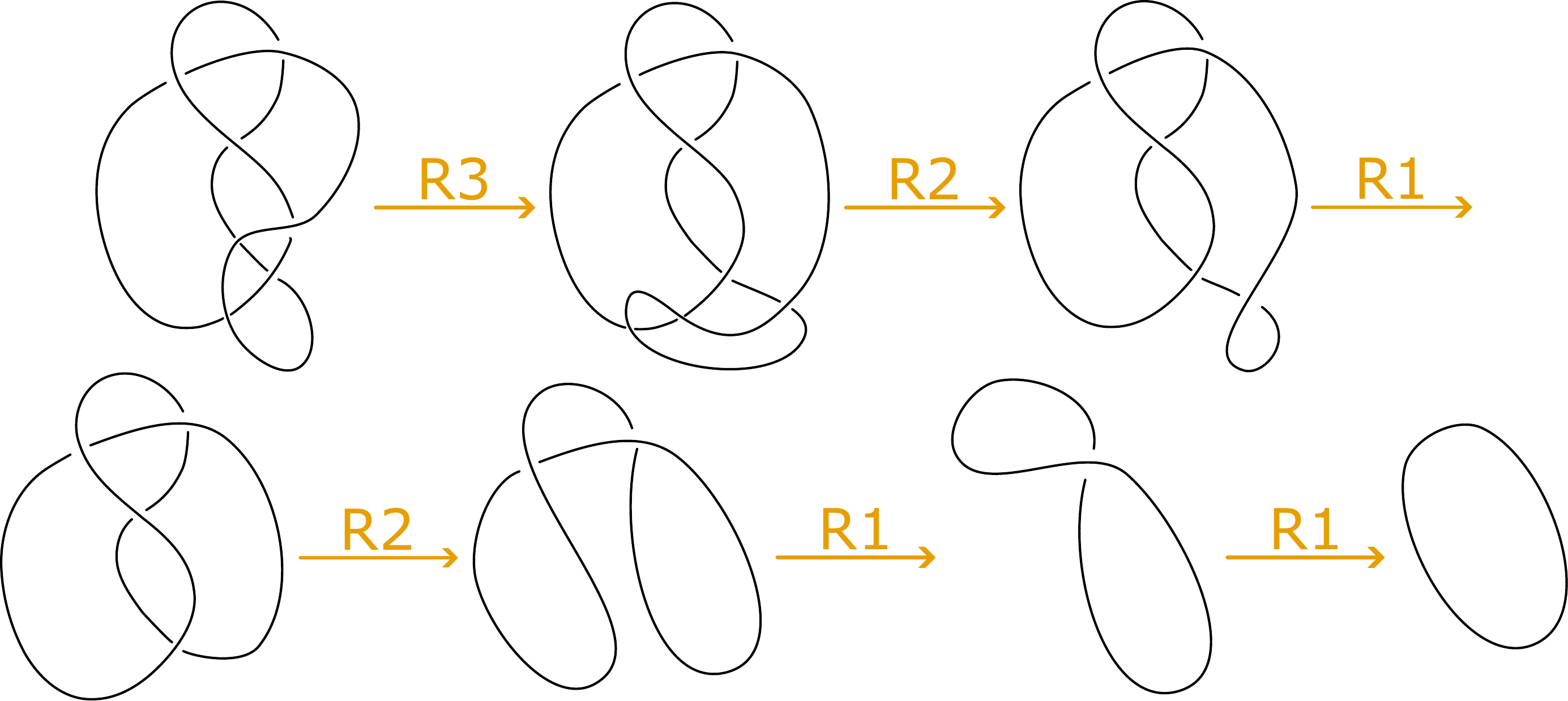}
            \caption{Deformaci\'on de un nudo al nudo trivial}
            \label{ej_redemeister}
        \end{figure}

Observemos el ejemplo de la Figura \ref{ej_redemeister}, si consideremos el nudo a izquierda y aplicamos la siguiente secuencia de movimientos de Reidemeister, obtenemos el nudo trivial. Este ejemplo nos permite preguntarnos si estos movimientos permiten determinar la equivalencia entre nudos, en efecto se obtiene el siguiente teorema. 

\teorema{Reidemeister}{
    Dos nudos $K_1,K_2$ son equivalentes si y solo si existe un diagrama $D_i$ con $i\in\{1,2\}$ para cada $K_i$ tales que $K_1$ se puede deformar en $K_2$ via una secuencia de movimientos de Redemeister.
}


Hasta el momento no se ha tenido en cuenta una orientaci\'on para nuestros nudos. En caso de que se requiera escoger una orientaci\'on para nuestro nudos, obtendremos los nudos denotados \textit{nudos orientados}, esta orientaci\'on se escoge en un punto inicial
y recorrer el nudo en una sola direcci\'on hasta regresar al punto inicial.\\

\begin{figure}[H]
    \centering
    \includegraphics[width=0.5\linewidth]{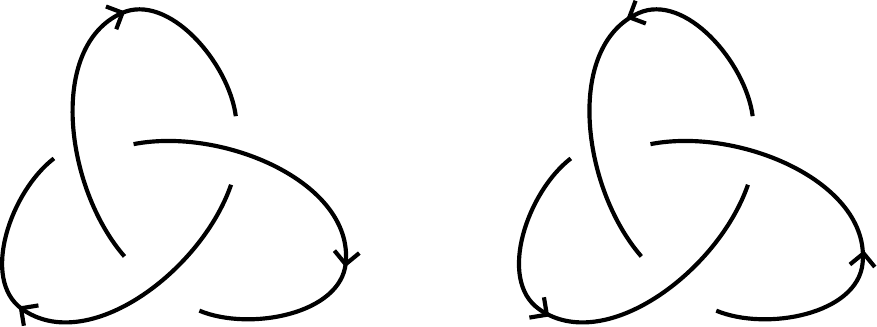}
    \caption{Orientaciones del nudo tr\'ebol}
    \label{trebol orientado}
\end{figure}

Por ejemplo, para el nudo tr\'ebol podemos considerar orientarlo en el sentido de las manecillas del reloj o en sentido contrario (ver Figura \ref{trebol orientado}). 

Esta noci\'on de orientaci\'on de un nudo nos servir\'a m\'as adelante para la construcci\'on de superficies asociadas al nudo.

\section{Familias de Nudos}\label{familias}
Los nudos se clasifican en tres familias disjuntas: t\'oricos, hiperb\'olicos y satelitales. En \'esta secci\'on introducimos de manera breve \'estas tres familias y describimos ejemplos de nudos que pertenecen a cada una de ellas. 

\subsection{Nudos T\'oricos}\label{toricos}
\begin{figure}[h!]
    \centering
    \includegraphics[width=0.4\linewidth]{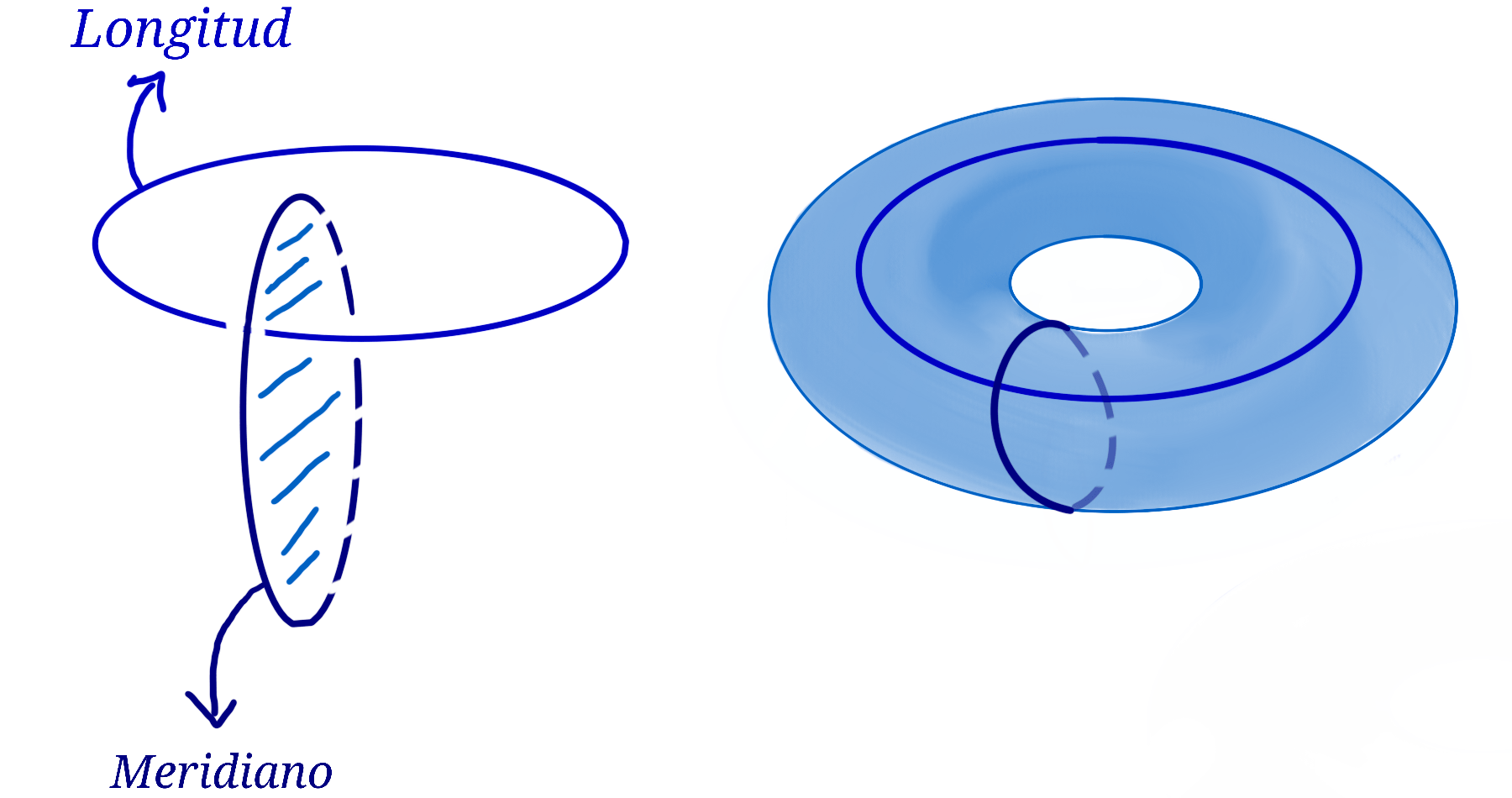}
    \caption{Toro a partir de meridianos y longitudinales}
    \label{toro s1 xs1}
\end{figure}

Los nudos t\'oricos son un tipo de nudos especiales que se pueden formar sobre la superficie de un toro $S^1\times S^1$ no anudado en $\mathbb{R}^3$. Para esta construcci\'on se debe tener en cuenta que el toro est\'a formado por dos tipos de c\'irculos, uno en direcci\'on meridional (el camino `corto' del toro) y otro c\'irculo en direcci\'on longitudinal (el camino que encierra el agujero del toro) como lo ilustra la \Cref{toro s1 xs1}.

\begin{figure}[h!]
    \centering
    \includegraphics[height=2cm]{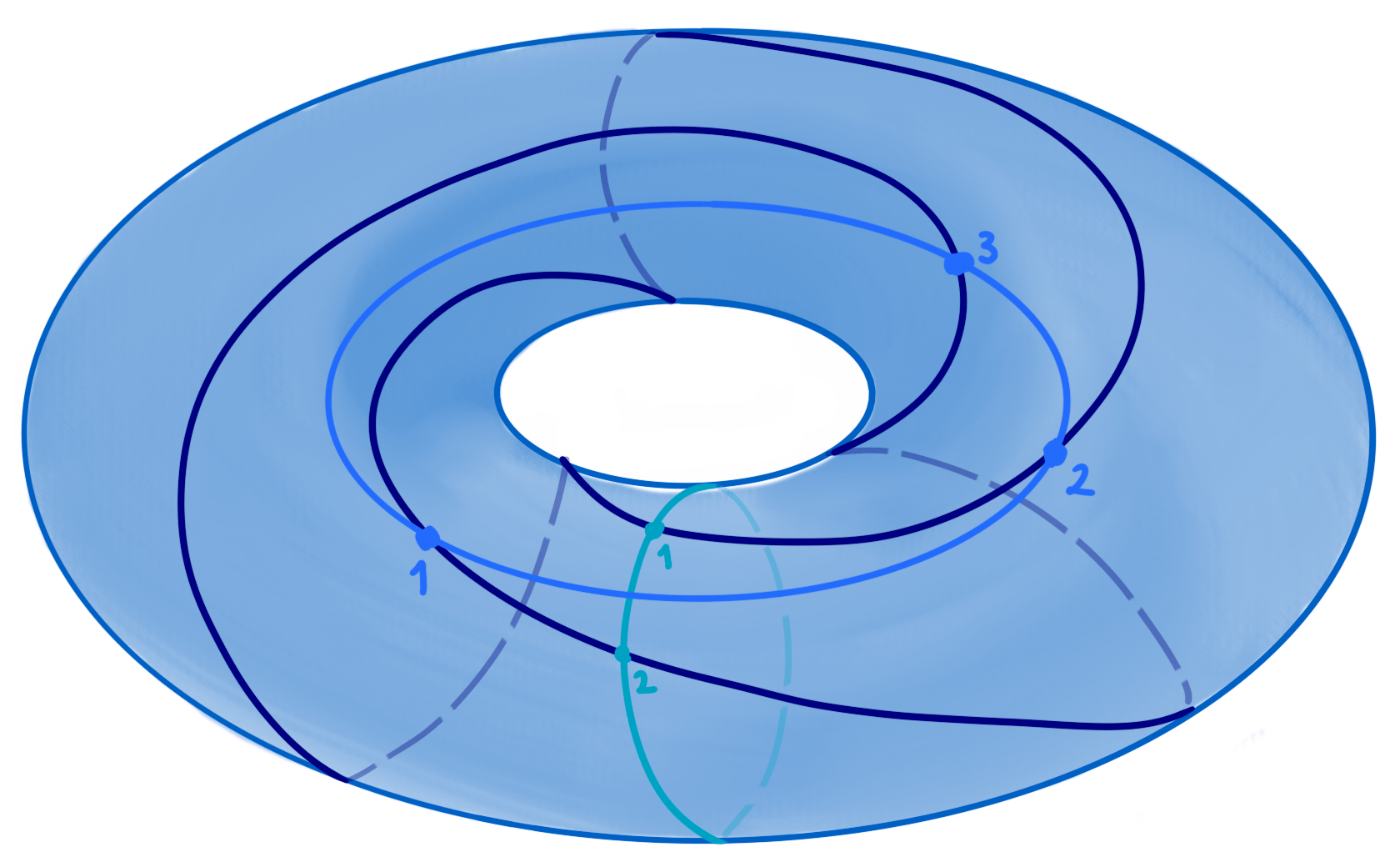}
    \caption{Nudo $T_{3,2}$ (Tr\'ebol)} 
    \label{trebol}
\end{figure}

Cada nudo t\'orico est\'a determinado por un par de enteros coprimos $p,q$ y se define como
\begin{equation}\label{T(p,q)}
T_{p,q}=\{(z^p,z^q) \mid z\in S^1\},
\end{equation}
es decir, es la curva cerrada que rodea $p$ veces el disco interno del toro s\'olido, y $q$ veces el agujero. La direcci\'on en la que el hilo del nudo se dispone en el toro se hace de tal forma que cuando $pq>0$, los hilos formen una h\'elice con torsi\'on a la derecha
En ocasiones es un poco dif\'icil determinar de manera directa si esta asignaci\'on de \textit{envolturas} es correcta, una forma de verificarlo es adicionando un meridiano y una curva longitudinal y observar cu\'antas veces estas cortan al nudo, indicando los valores de $q$ y $p$ respectivamente. Rec\'iprocamente si el nudo cruza $p$ veces la curva longitudinal, entonces el nudo debe envolverse en el sentido meridional $p$ veces. Por ejemplo para el nudo tr\'ebol, en el toro este corresponde al nudo t\'orico $T_{3,2}$ (ver Figura \ref{trebol}).\\

\begin{figure}[H]
    \centering
    \includegraphics[height=2cm]{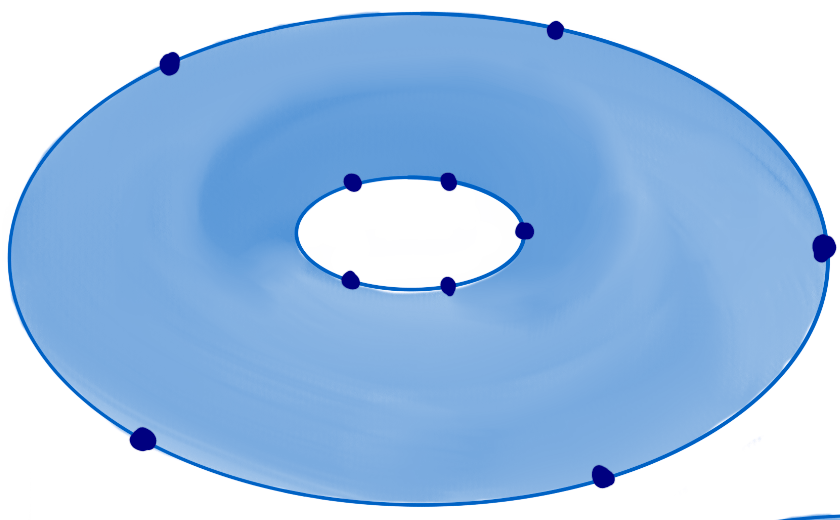}\hfill%
    \includegraphics[height=2cm]{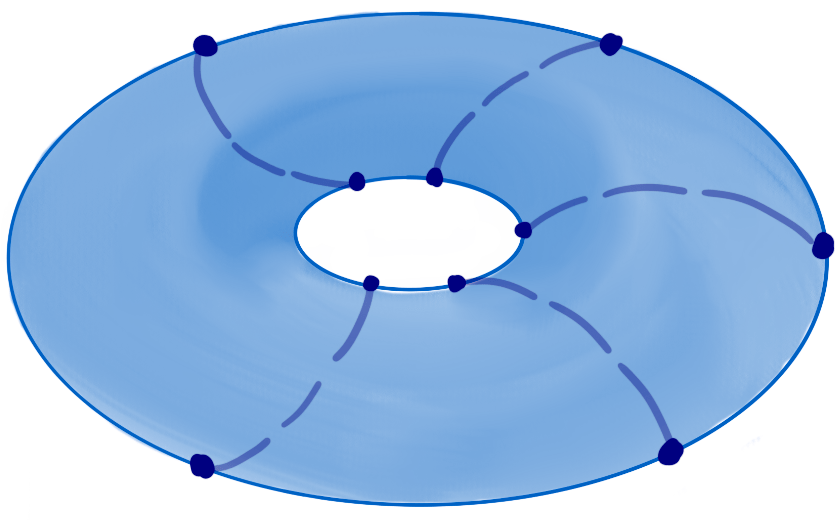}\hfill%
    \includegraphics[height=2cm]{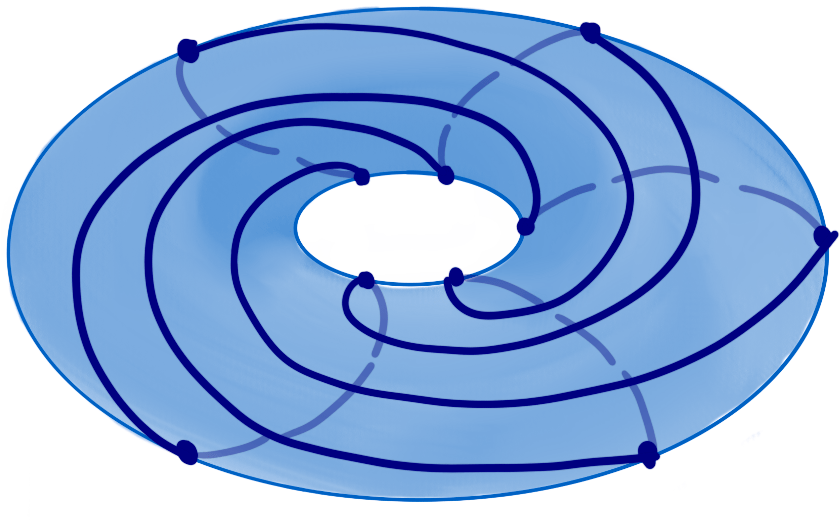}
    \caption{Construcci\'on del toro $T_{5,3}$}
    \label{construccion en el toro}
\end{figure}

Los nudos t\'oricos se pueden construir a partir de un peque\~no algoritmo. Para construir el nudo t\'orico $T_{p,q}$, se marcan $|p|$ puntos tanto en la curva longitudinal exterior como en la interior, y se unen con un arco los puntos directamente opuestos por la parte inferior del toro. En la parte superior del toro se une cada punto exterior con el punto interior que est\'a a $q/p$ de distancia en el sentido de las manecillas del reloj, es decir saltamos $q$ puntos y conectamos en el $q+1$. La \Cref{construccion en el toro} ilustra la construcci\'on del nudo t\'orico $T_{5,3}$. \\ 

\begin{figure}[h!]
    \centering
    \includegraphics[width=0.3\linewidth]{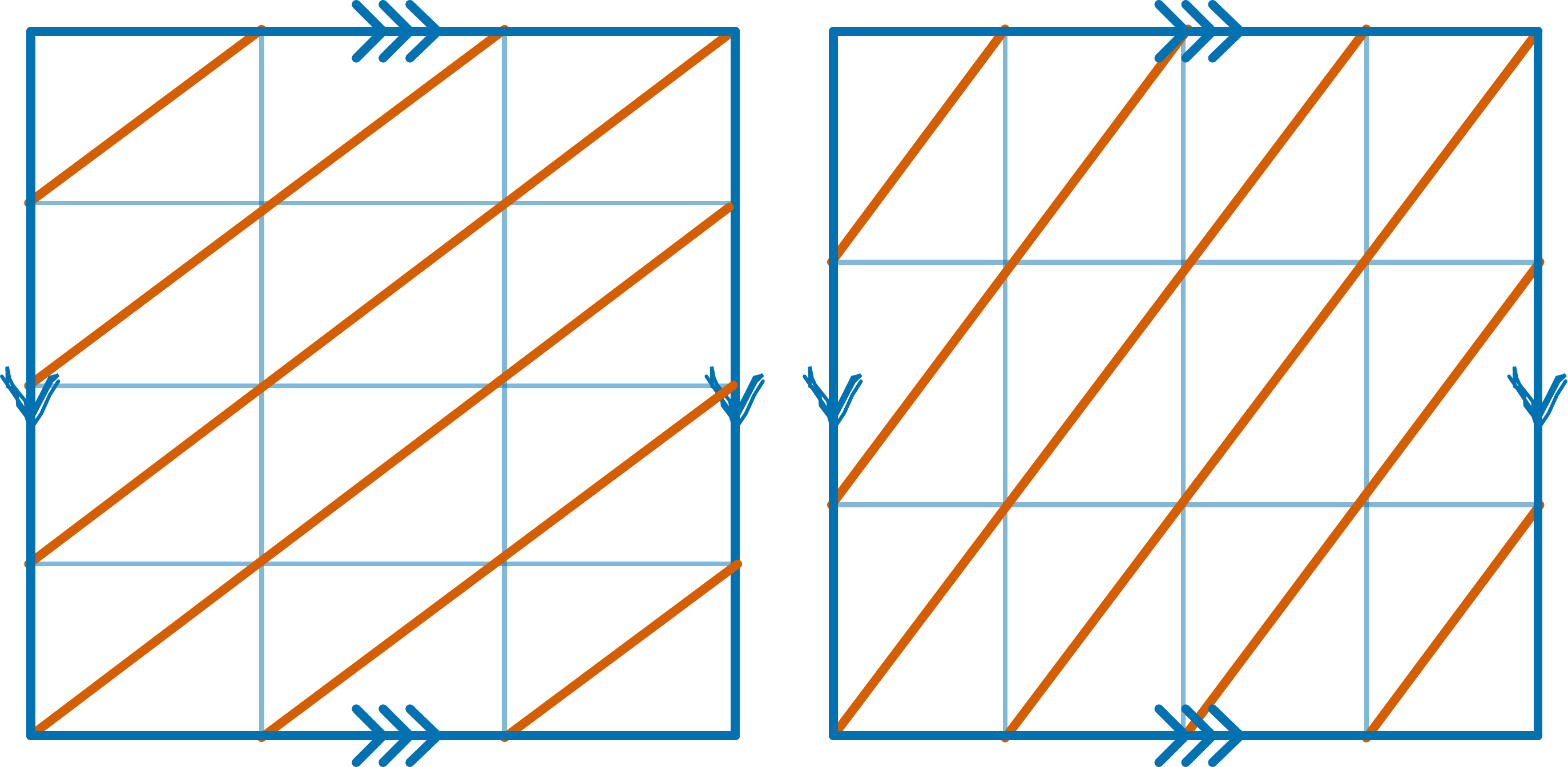}
    \caption{Diagrama poligonal del $T_{3,4}$ y $T_{4,3}$}
    \label{diagrama en toro}
\end{figure}

Para representar los nudos t\'oricos en el diagrama poligonal can\'onico del toro (ver \Cref{identificacion_toro}), primero se dibuja una cuadr\'icula de tama\~no $|p|\times |q|$ en el cuadrado. Si $pq>0$ se traza una l\'inea recta que une la esquina inferior izquierda con la esquina superior derecha de cada rect\'angulo de la cuadr\'icula. Si $pq<0$ se unen la esquina superior izquierda con la inferior derecha. La \Cref{diagrama en toro} muestra la representaci\'on de los nudos $T_{3,4}$ y $T_{4,3}$. En general, tenemos la siguiente proposici\'on.

\proposicion{}{Sean $p,q$ primos relativos. Los nudos $T_{p,q}$ y $T_{q,p}$ son equivalentes. \label{isotopia toricos
}}

\begin{figure}[h!]
    \centering
    \includegraphics[width=0.45\linewidth]{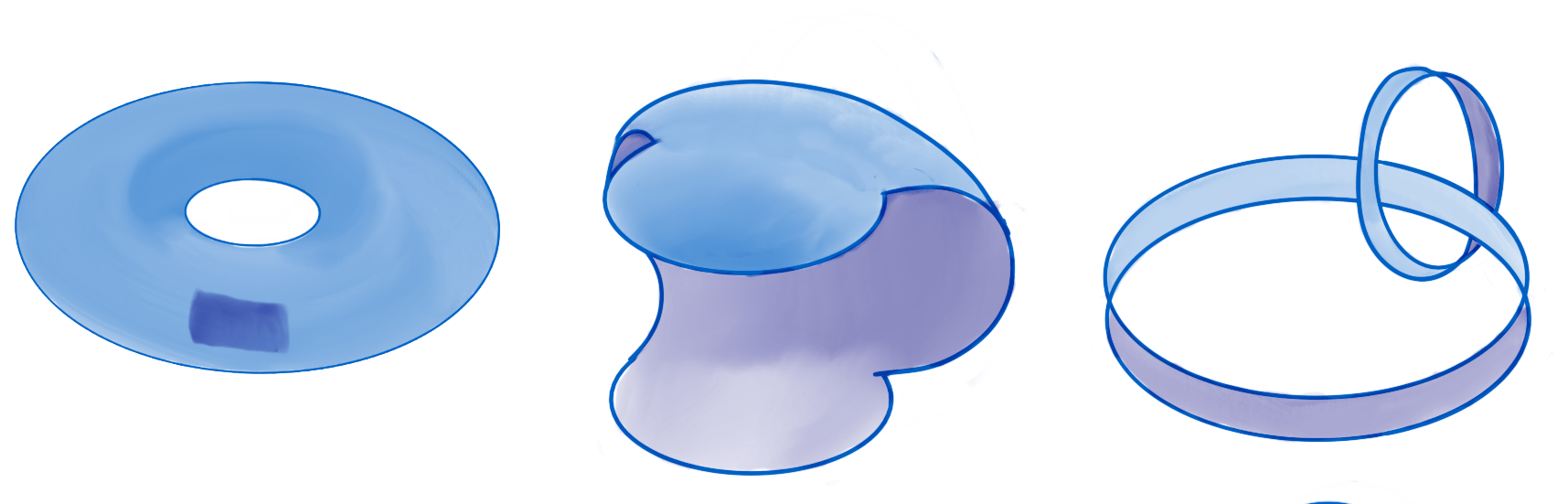}\hspace{1ex}%
    \includegraphics[width=0.45\linewidth]{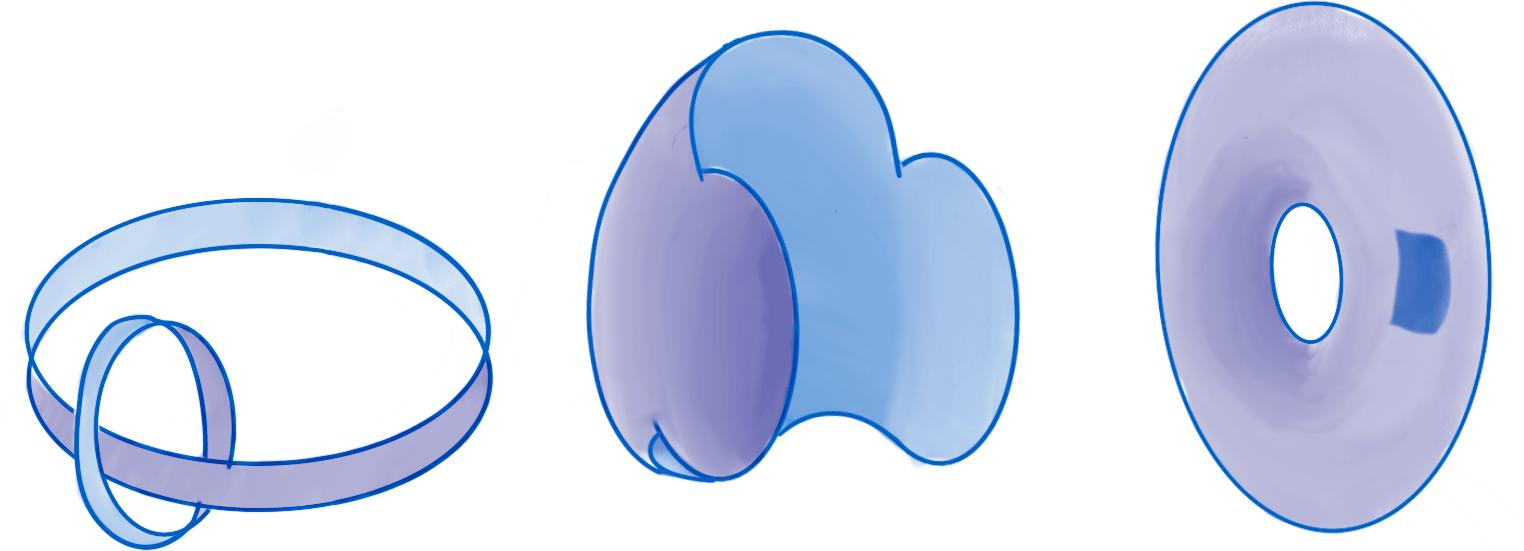}
    \caption{Inversi\'on del toro, cambiando curvas longitudinales y meridianos de posici\'on}
    \label{inversion toro}
\end{figure}

Para ver la isotop\'ia eliminamos un disco punteado de la superficie del toro, asegur\'andonos que el disco no toca el nudo, luego podemos observar que al \textit{agrandar} el hueco generado por el disco eliminado, el toro se convierte en dos bandas, una relacionada con el meridiano y una con la longitudinal, posteriormente las bandas cambian de roles a partir de enviar hacia adentro el meridiano y hacia afuera la longitudinal, logrando el toro con las bandas invertidas, as\'i obteniendo el nudo t\'orico $T_{q,p}$. (ver \Cref{inversion toro}).\\

\subsection{Nudos Hiperb\'olicos}\label{hiperbolicos}
\begin{figure}[h]
\centering
\includegraphics[height=2.5cm]{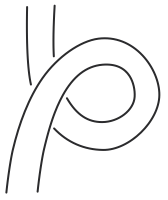}\hspace{1cm}
\includegraphics[height=2.5cm]{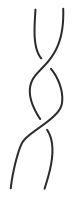}\hspace{1cm}
\includegraphics[height=2.5cm]{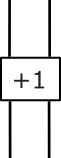}
\caption{Tres representaciones equivalentes de bandas con giros completos. }\label{bandas}
\end{figure}

\begin{figure}[h]
\centering
\def\svgwidth{0.2\textwidth}
\begingroup%
  \makeatletter%
  \providecommand\color[2][]{%
    \errmessage{(Inkscape) Color is used for the text in Inkscape, but the package 'color.sty' is not loaded}%
    \renewcommand\color[2][]{}%
  }%
  \providecommand\transparent[1]{%
    \errmessage{(Inkscape) Transparency is used (non-zero) for the text in Inkscape, but the package 'transparent.sty' is not loaded}%
    \renewcommand\transparent[1]{}%
  }%
  \providecommand\rotatebox[2]{#2}%
  \newcommand*\fsize{\dimexpr\f@size pt\relax}%
  \newcommand*\lineheight[1]{\fontsize{\fsize}{#1\fsize}\selectfont}%
  \ifx\svgwidth\undefined%
    \setlength{\unitlength}{134.10589023bp}%
    \ifx\svgscale\undefined%
      \relax%
    \else%
      \setlength{\unitlength}{\unitlength * \real{\svgscale}}%
    \fi%
  \else%
    \setlength{\unitlength}{\svgwidth}%
  \fi%
  \global\let\svgwidth\undefined%
  \global\let\svgscale\undefined%
  \makeatother%
  \begin{picture}(1,2.38852268)%
    \lineheight{1}%
    \setlength\tabcolsep{0pt}%
    \put(0,0){\includegraphics[width=\unitlength,page=1]{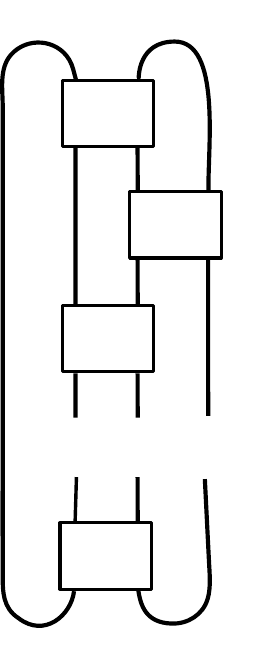}}%
    \put(0.274205,1.93231102){\color[rgb]{0.12941176,0.12941176,0.12941176}\makebox(0,0)[lt]{\lineheight{1.25}\smash{\begin{tabular}[t]{l}$a_1$\end{tabular}}}}%
    \put(0.52773151,1.54677497){\color[rgb]{0.12941176,0.12941176,0.12941176}\makebox(0,0)[lt]{\lineheight{1.25}\smash{\begin{tabular}[t]{l}$-a_2$\end{tabular}}}}%
    \put(0.274205,1.12530723){\color[rgb]{0.12941176,0.12941176,0.12941176}\makebox(0,0)[lt]{\lineheight{1.25}\smash{\begin{tabular}[t]{l}$a_3$\end{tabular}}}}%
    \put(0.29845415,0.34436118){\color[rgb]{0.12941176,0.12941176,0.12941176}\makebox(0,0)[lt]{\lineheight{1.25}\smash{\begin{tabular}[t]{l}$a_l$\end{tabular}}}}%
    \put(0.47678804,0.75606573){\color[rgb]{0.12941176,0.12941176,0.12941176}\makebox(0,0)[lt]{\lineheight{1.25}\smash{\begin{tabular}[t]{l}$\vdots$\end{tabular}}}}%
  \end{picture}%
\endgroup%

\def\svgwidth{0.2\textwidth}
\begingroup%
  \makeatletter%
  \providecommand\color[2][]{%
    \errmessage{(Inkscape) Color is used for the text in Inkscape, but the package 'color.sty' is not loaded}%
    \renewcommand\color[2][]{}%
  }%
  \providecommand\transparent[1]{%
    \errmessage{(Inkscape) Transparency is used (non-zero) for the text in Inkscape, but the package 'transparent.sty' is not loaded}%
    \renewcommand\transparent[1]{}%
  }%
  \providecommand\rotatebox[2]{#2}%
  \newcommand*\fsize{\dimexpr\f@size pt\relax}%
  \newcommand*\lineheight[1]{\fontsize{\fsize}{#1\fsize}\selectfont}%
  \ifx\svgwidth\undefined%
    \setlength{\unitlength}{146.51204988bp}%
    \ifx\svgscale\undefined%
      \relax%
    \else%
      \setlength{\unitlength}{\unitlength * \real{\svgscale}}%
    \fi%
  \else%
    \setlength{\unitlength}{\svgwidth}%
  \fi%
  \global\let\svgwidth\undefined%
  \global\let\svgscale\undefined%
  \makeatother%
  \begin{picture}(1,2.18838833)%
    \lineheight{1}%
    \setlength\tabcolsep{0pt}%
    \put(0,0){\includegraphics[width=\unitlength,page=1]{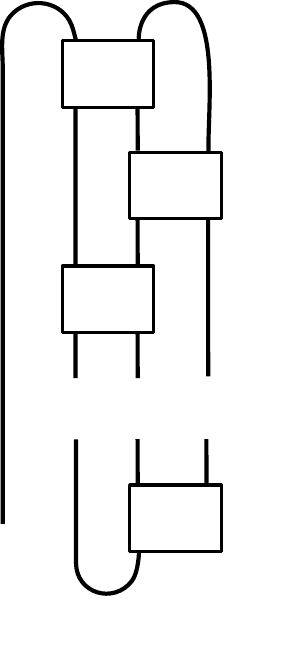}}%
    \put(0.25098622,1.90008273){\color[rgb]{0.12941176,0.12941176,0.12941176}\makebox(0,0)[lt]{\lineheight{1.25}\smash{\begin{tabular}[t]{l}$a_1$\end{tabular}}}}%
    \put(0.48304494,1.54719261){\color[rgb]{0.12941176,0.12941176,0.12941176}\makebox(0,0)[lt]{\lineheight{1.25}\smash{\begin{tabular}[t]{l}$-a_2$\end{tabular}}}}%
    \put(0.25098622,1.16141337){\color[rgb]{0.12941176,0.12941176,0.12941176}\makebox(0,0)[lt]{\lineheight{1.25}\smash{\begin{tabular}[t]{l}$a_3$\end{tabular}}}}%
    \put(0.45614274,0.45193728){\color[rgb]{0.12941176,0.12941176,0.12941176}\makebox(0,0)[lt]{\lineheight{1.25}\smash{\begin{tabular}[t]{l}$-a_{l}$\end{tabular}}}}%
    \put(0,0){\includegraphics[width=\unitlength,page=2]{2-puente-even.pdf}}%
    \put(0.45653191,0.81879247){\color[rgb]{0.12941176,0.12941176,0.12941176}\makebox(0,0)[lt]{\lineheight{1.25}\smash{\begin{tabular}[t]{l}$\vdots$\end{tabular}}}}%
  \end{picture}%
\endgroup%

\caption{Nudos de 2 puentes con $l$ parametros. En la izquierda se ilustra el caso $l$ impar, a la derecha el caso $l$ par.}\label{2-bridge}
\end{figure}

Un nudo $K$ es hiperb\'olico si su complemento $S^3\setminus K$ admite una m\'etrica riemanniana de curvatura negativa constante. La clase de nudos de 2 puentes que no son toroidales es una familia infinita de nudos hiperb\'olicos \cite{menasco}. Definiremos esta familia de nudos de manera diagram\'atica. Dados numeros enteros $a_1, a_2,\ldots,a_l$ definimos el nudo $K(a_1, a_2,\ldots,a_l)$ como en la \Cref{2-bridge}, dependiendo de la paridad de $l$. Los cuadrados representan el n\'umero de giros completos que se le dan a las hebras paralelas. La \Cref{bandas} ilustra el caso de giros positivos. Los giros negativos se obtienen cambiando hebras superiores por hebras inferiores en las im\'agenes que aparecen en la derecha y el centro de la \Cref{bandas}.

\dfe{Nudos de 2 puentes}{Un nudo $K$ es un nudo de dos puentes si admite un diagrama como en la \Cref{2-bridge}.}

Existe una biyecci\'on entre el conjunto de nudos de 2 puentes y el conjunto de n\'umeros racionales con numerador impar. El n\'umero racional asociado a un nudo de 2 puentes se denomina forma normal de Schubert \cite{schubert-2bridge}. De manera expl\'icita, al nudo $K(a_1, a_2,\ldots,a_l)$ se asocia el n\'umero racional con fracci\'on continua $$[a_1, a_2,\ldots,a_l]=a_1+\cfrac{1}{a_2+\cfrac{1}{\ddots \,\cfrac{1}{a_{l-1}+\cfrac{1}{a_l}}}}.$$

Dentro de la familia de nudos de 2 puentes con dos par\'ametros, destacamos los nudos $C_m=K(2m,2m)$ con fracci\'on $\tfrac{4m^2+1}{2m}$, $m\in\Z$. 

\proposicion{Nudos Anfiquirales}{\label{anfiquiral}
Para todo $m\in\Z$, los nudos $C_m$ y $C_{-m}$ son equivalentes.}

La demostraci\'on de la proposici\'on se la dejamos al lector interesado. 

\begin{figure}[h]
\centering
\includegraphics[height=2.75cm]{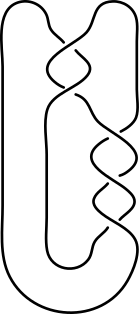}\qquad
\includegraphics[height=2.75cm]{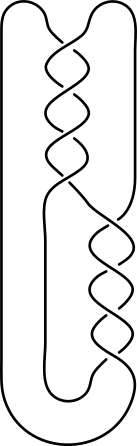}
\caption{Ejemplos de nudos de 2-puentes con dos parametros. A la izquierda se muestra el nudo $5_2=K(2,3)$ y a la derecha el nudo $8_3=K(4,4)$.}\label{double_twist}
\end{figure}

\subsection{Nudos Satelitales}\label{satelitales}
\begin{figure}[h]
\centering
\includegraphics[height=3cm]{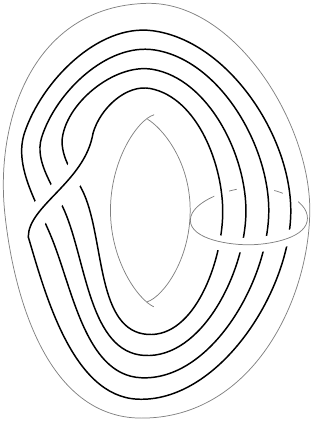}\qquad\qquad%
\includegraphics[height=3cm]{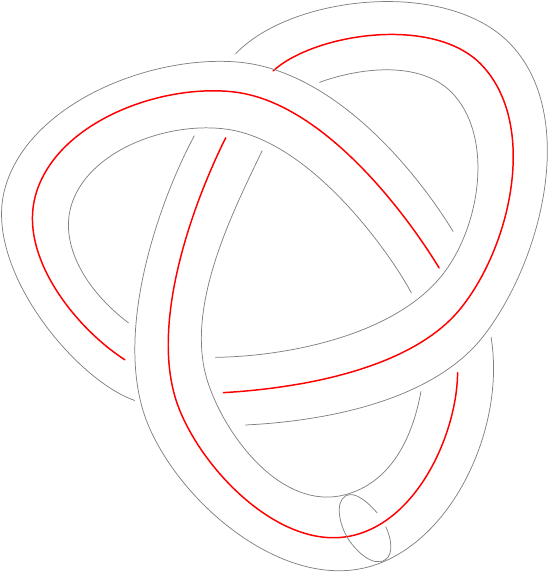}\qquad\qquad%
\includegraphics[height=3cm]{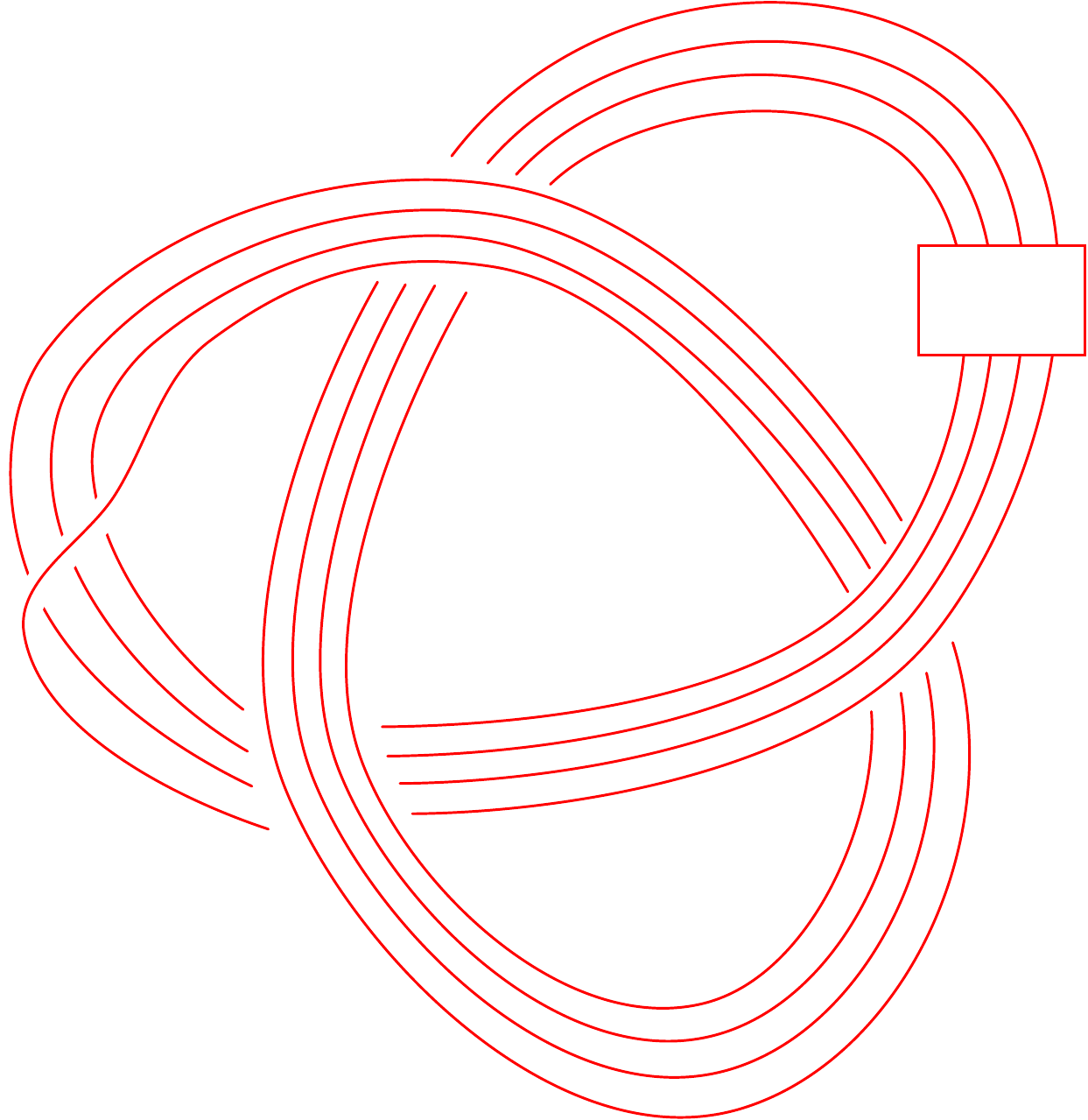}
\caption{Los nudos que hacen parte de una operaci\'on satelital. Patr\'on $P=C_{4,1}$ (izquierda), \'orbita $T_{2,3}$(centro) y nudo satelital $C_{4,1}\left(T_{2,3}\right)$ (derecha). El cuadro de la derecha contiene tres giros completos negativos que corrigen la torsi\'on (o writhe) $+3$ en el diagrama de $T_{2,3}$.}\label{satelite}
\end{figure}

En una operaci\'on satelital dos nudos (un patr\'on y una \'orbita) se combinan para producir un tercer nudo (el sat\'elite) como en \cref{satelite}. En t\'erminos de espacios exteriores o complementos, un nudo $J$ es un nudo sat\'elite o satelital si y solo si su espacio exterior $S^3\setminus N(J)$ contiene un toro incompresible que no es paralelo a la frontera. La clase de nudos satelitales incluye las sumas conexas de nudos. 

\begin{definition} [Operaciones Satelitales]
Sea $P \subset D^2\times S^1$ un nudo encajado en un toro s\'olido y $K \subset S^3$ un nudo. Si $h:D^2\times S^1\to S^3$ es un embebimiento que identifica el toro s\'olido $D^2\times S^1$ con una vecindad tubular de $K$, entonces $h(P)$ es el nudo satelital $P(K)$ con \'orbita $K$ y patr\'on $P$.\\
El n\'umero de \'indice $w$ de la operaci\'on es el n\'umero de intersecci\'on $\# (P\cap D^2\times \{*\})$. 
\end{definition}

Es natural preguntarse c\'omo la complejidad del nudo $P(K)$ est\'a determinada por las de sus piezas constituyentes $P$ y $K$. En t\'erminos del 3-g\'enero, Schubert~\cite{schubert} da una f\'ormula exacta:

\teorema{Schubert}{Para cualquier nudo no trivial $K$,
\[g_3(P(K))= g_3(P)+ |w| g_3(K),\]}

\section{Superficies con una componente conexa en su frontera}\label{superficies_frontera}
A lo largo de este cap\'itulo estudiaremos c\'omo podemos asociar una superficie a un nudo de forma que a partir de esta asignaci\'on podamos extraer informaci\'on sobre la equivalencia de nudos.

\subsection{Clasificaci\'on}
\begin{figure}[h!]
    \centering
    \includegraphics[width=0.5\linewidth]{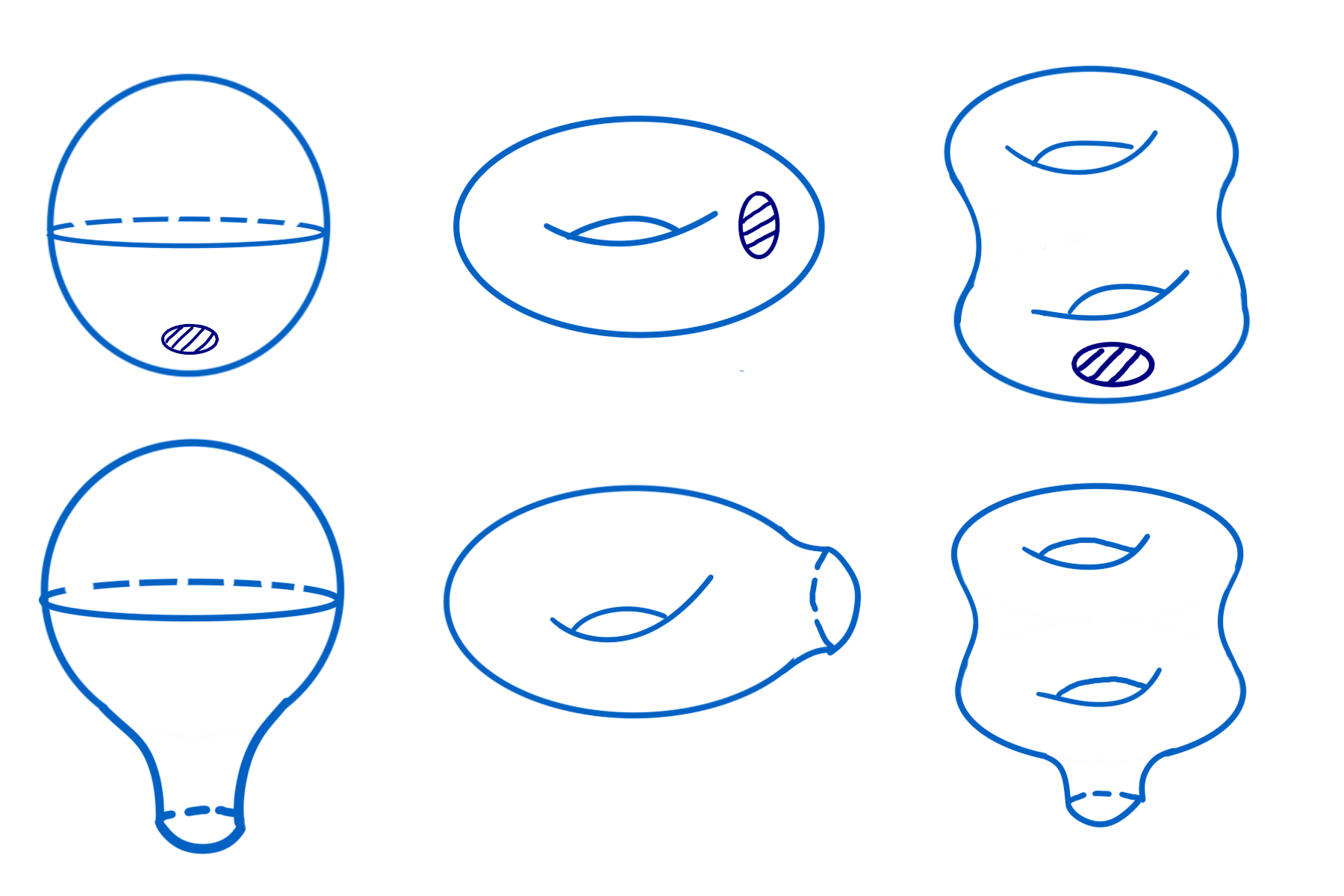}
    \caption{Superficies cerradas, conexas y con una componente conexa en su frontera}
    \label{frontera}
\end{figure}

En las secciones anteriores estuvimos interesados en estudiar las superficies cerradas y conexas, ahora queremos hacer una peque\~na modificaci\'on y permitir la presencia de una \'unica componente conexa en la frontera. Esto se puede lograr f\'acilmente al remover un disco de cada una de las superficies del Teorema \ref{homeomorfismo}, obteniendo superficies como se muestran en la Figura \ref{frontera}. Asimismo, una superficie $F$ compacta, orientable y con frontera conexa da a lugar a una superficie cerrada y orientable de la siguiente manera: sea $g:\partial F\to S^1=\partial D^2$ la funci\'on que identifica la componente conexa de la frontera de $F$ con el c\'irculo $S^1$. Se forma entonces $\widehat{F}$ como el espacio cociente de $F\sqcup D^2$ con relaci\'on de equivalencia definida por $x\sim g(x)$ para todo $x\in \partial F$ (y por supuesto todo $g(x)\in S^1=\partial D^2$). En otras palabras, hay una correspondencia uno-a-uno entre las superficies orientables cerradas, y las superficies orientables con frontera conexa. De \'esta correspondencia se sigue el siguiente resultado.

\teorema{Homeomorfismo entre superficies}{\label{clas-frontera}
Sean $F_1$ y $F_2$ dos superficies compactas, orientables y con una componente conexa en su frontera. $F_1$ y $F_2$ son \textbf{homeomorfas} si y solamente si tienen la misma caracter\'istica de Euler.}

En t\'erminos de la caracter\'istica de Euler tenemos $\chi(\widehat{F})=\chi(F)+\chi(D^2)-\chi(S^1)=\chi(F)+1$. Por lo tanto \cref{euler-g} muestra que 
\begin{equation}\label{euler-g-frontera}
\chi(F)=1-2g(F).
\end{equation}

\subsection{Superficies de Seifert}

En esta secci\'on queremos relacionar las superficies compactas, orientables y con una componente conexa en su frontera, con los nudos de tal forma que el nudo sea justamente la frontera de alguna superficie. El siguiente teorema nos garantiza la existencia de dicha superficie para un nudo contenido en $\mathbb{R}^3$.

\teorema{Seifert 1926}{Dado un diagrama $D$ que representa un nudo $K\subset\mathbb{R}^3$, existe una superficie $F\subset\mathbb{R}^3$ encajada, compacta y orientable tal que  $\partial F=K$. Dicha superficie se llama Superficie de Seifert del nudo $K$.}

La existencia de superficies de Seifert se puede establecer de varias maneras. Una de ellas es el llamado `algoritmo de Seifert', un proceso que construye una superficie a partir de un diagrama dado.\\

\begin{figure}[H]
    \centering
    \includegraphics[clip, trim=20 0 20 0, width=0.3\linewidth]{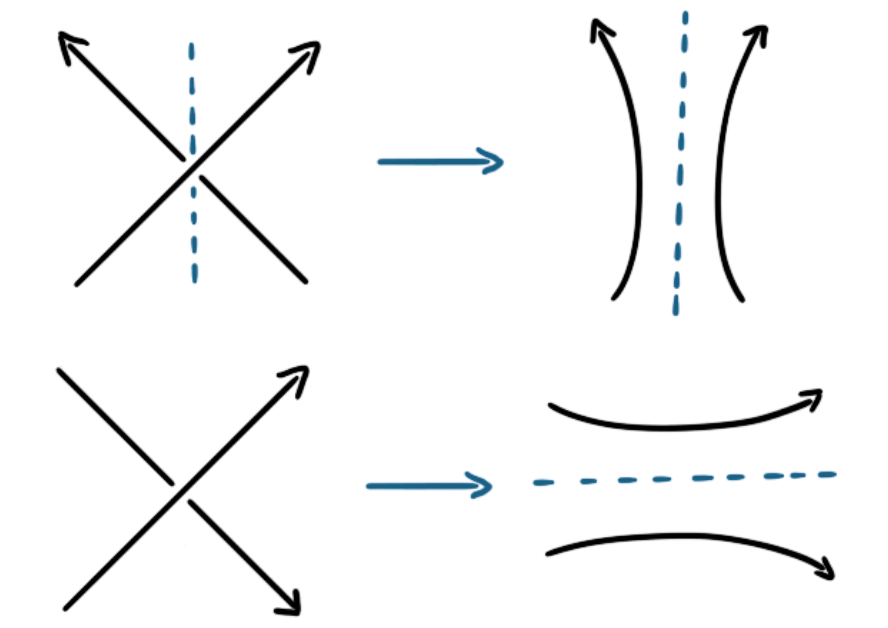}
    \caption{Resoluci\'on orientada de un cruce.}
    \label{res-cruce}
\end{figure}

\texttt{\noindent Algoritmo de Seifert.\\
Input: Un diagrama $D$ para un nudo $K$.
\begin{enumerate}\label{seifert}
\item Asignar una orientaci\'on al nudo.
\item Reemplazar cada cruce con su resoluci\'on orientada, como se muestra en la \Cref{res-cruce}. El diagrama resultante es una colecci\'on de c\'irculos orientados, conocidos como c\'irculos de Seifert. 
A cada c\'irculo de Seifert $C$ se le asigna la \textit{`altura de anidaci\'on'} dentada por $h(C)$ y definida como el n\'umero de otros c\'irculos que uno debe cruzar para llegar desde $C$ al punto `en el infinito' en $\R^2$.
\item Rellenar cada c\'irculo de Seifert $C$ con un disco en $\R^3$ paralelo al plano $xy$ y situado a la altura $z=h(C)$. Retomando la idea de orientaciones, pensamos en los discos como echos en tela con un lado oscuro y uno claro. Si $C$ est\'a orientado en el sentido de las manecillas del reloj, el lado visible del disco es el `oscuro', de lo contrario el lado visible es el `claro'. 
\item Conectar los discos usando bandas con giros en el lugar correspondiente a cada cruce. La orientaci\'on del giro de la banda se deriva del tipo de cruce. Como en el paso anterior, tomamos bandas echas de la misma tela de los discos. Siendo as\'i el caso, al darle media vuelta a la banda rectangular vemos un color distinto a cada lado. Usando la selecci\'on de discos del paso anterior nos podemos asegurar que en cada cruce los colores de la banda sean compatibles con los colores de los discos a los que se sujeta.
\end{enumerate}
Output: Una superficie orientada $S$ con frontera el nudo $K$.
}

Ahora, vamos a ejemplificar los pasos que se usan en el algoritmo de Seifert para obtener una superfie tal que la frontera de esta sea un nudo dado. 

\begin{figure}[H]
    \centering
    \includegraphics[clip, trim=50 500 50 400, width=0.5\linewidth]{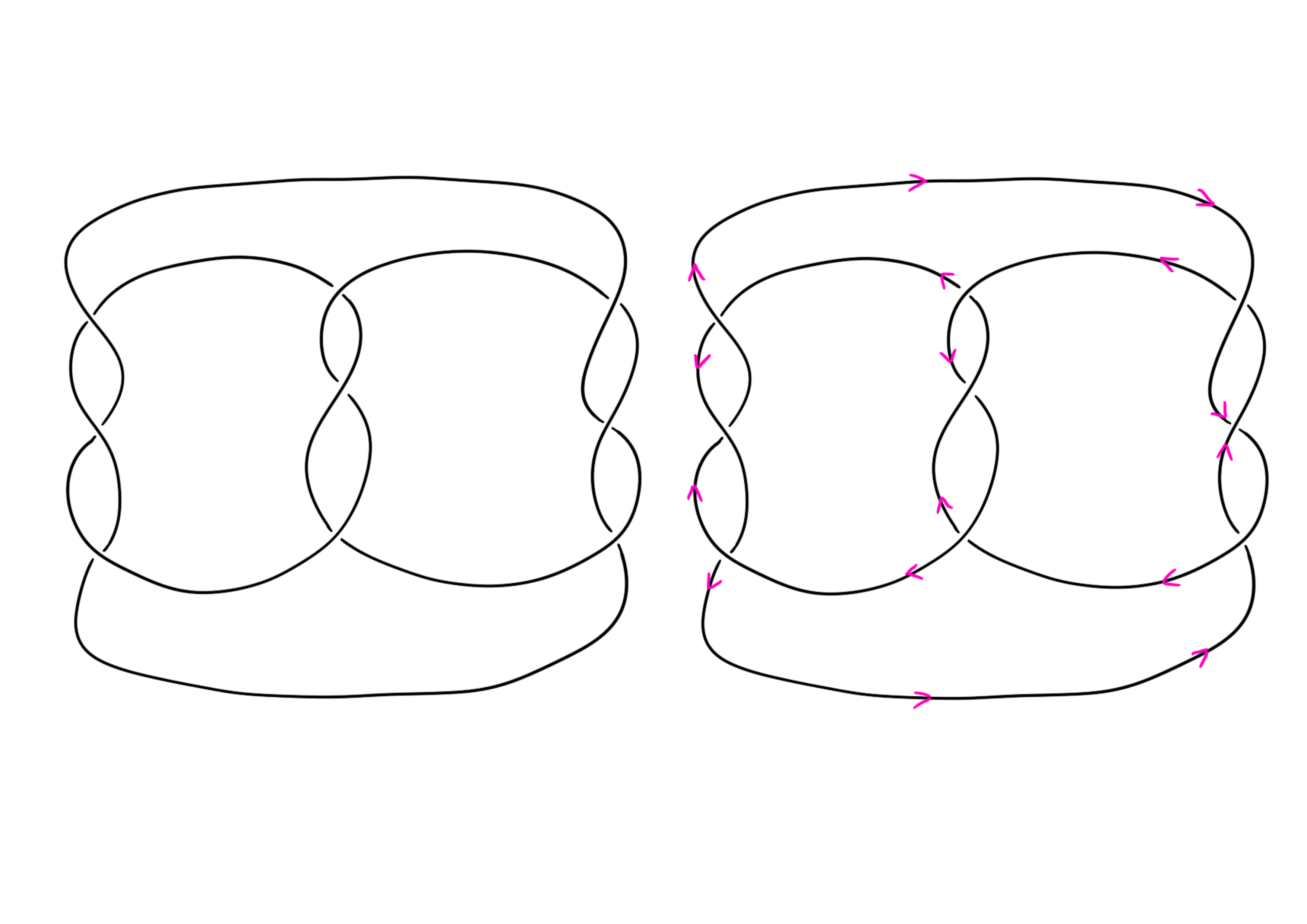}
    \caption{Orientaci\'on de un nudo}
    \label{seifertsup1}
\end{figure}

La Figura \ref{seifertsup1} muestra el resultado del primer paso, asociarle una orientaci\'on a nuestro nudo. En este caso, decidimos empezar por la derecha. Para efectuar el paso 2, dibuj\'amos un hilo adyacente al diagrama del nudo, y al llegar a un cruce lo resolvemos de manera que se preserve la orientaci\'on asignada. La forma de resolver los cortes se puede ver en la Figura \ref{res-cruce}.

\begin{figure}[H]
    \centering
    \includegraphics[clip, trim=50 500 50 450, width=0.5\linewidth]{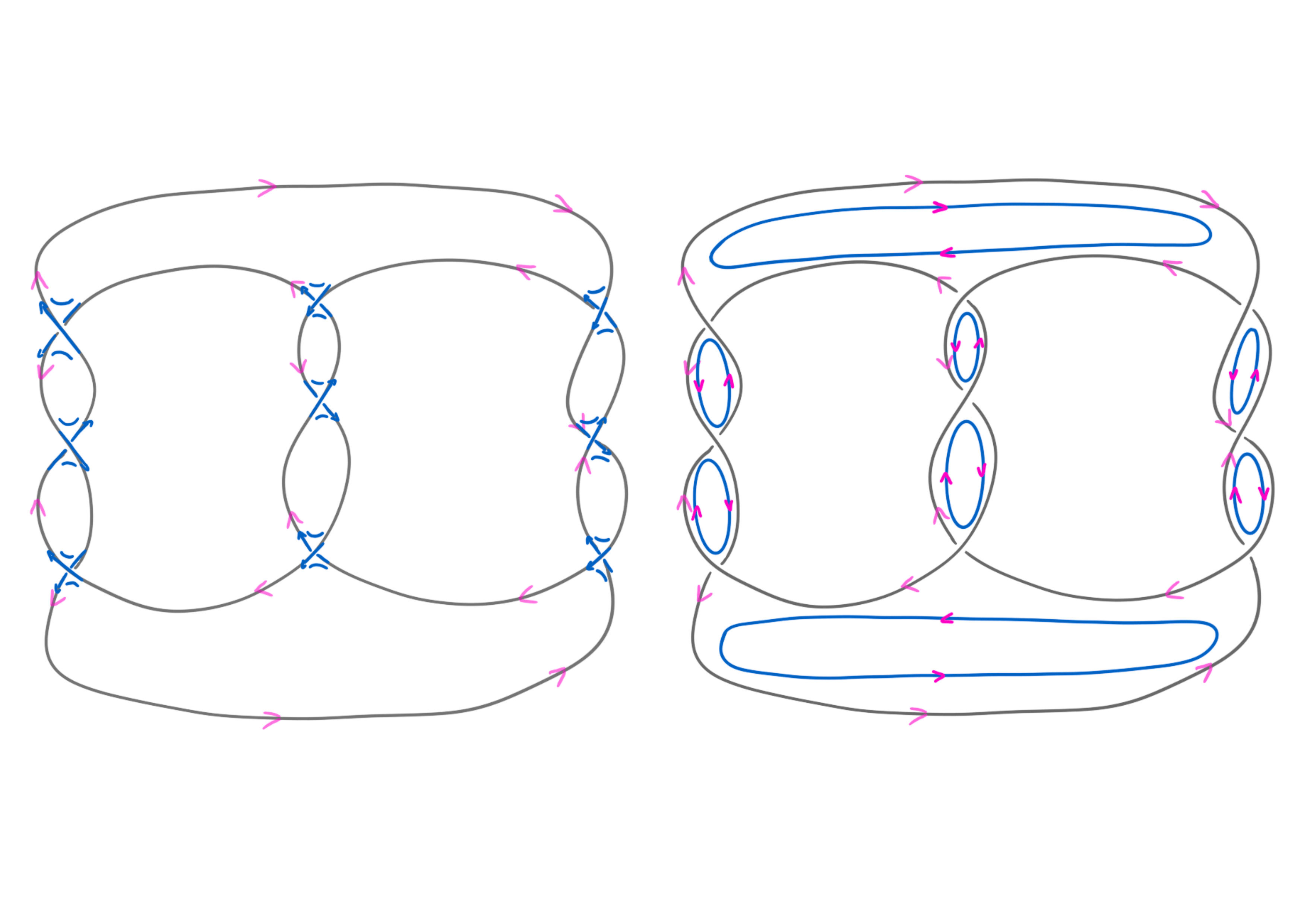}
    \caption{C\'irculos orientados}
    \label{seifertsup2}
\end{figure}

El resultado de estos hilos adyacentes a mi diagrama es una colecci\'on de c\'irculos orientados en $\mathbb{R}^2$, generalmente conocidos como \textit{ciclos de Seifert.}

\begin{figure}[H]
    \centering
    \includegraphics[clip, trim=50 500 50 450, width=0.6\linewidth]{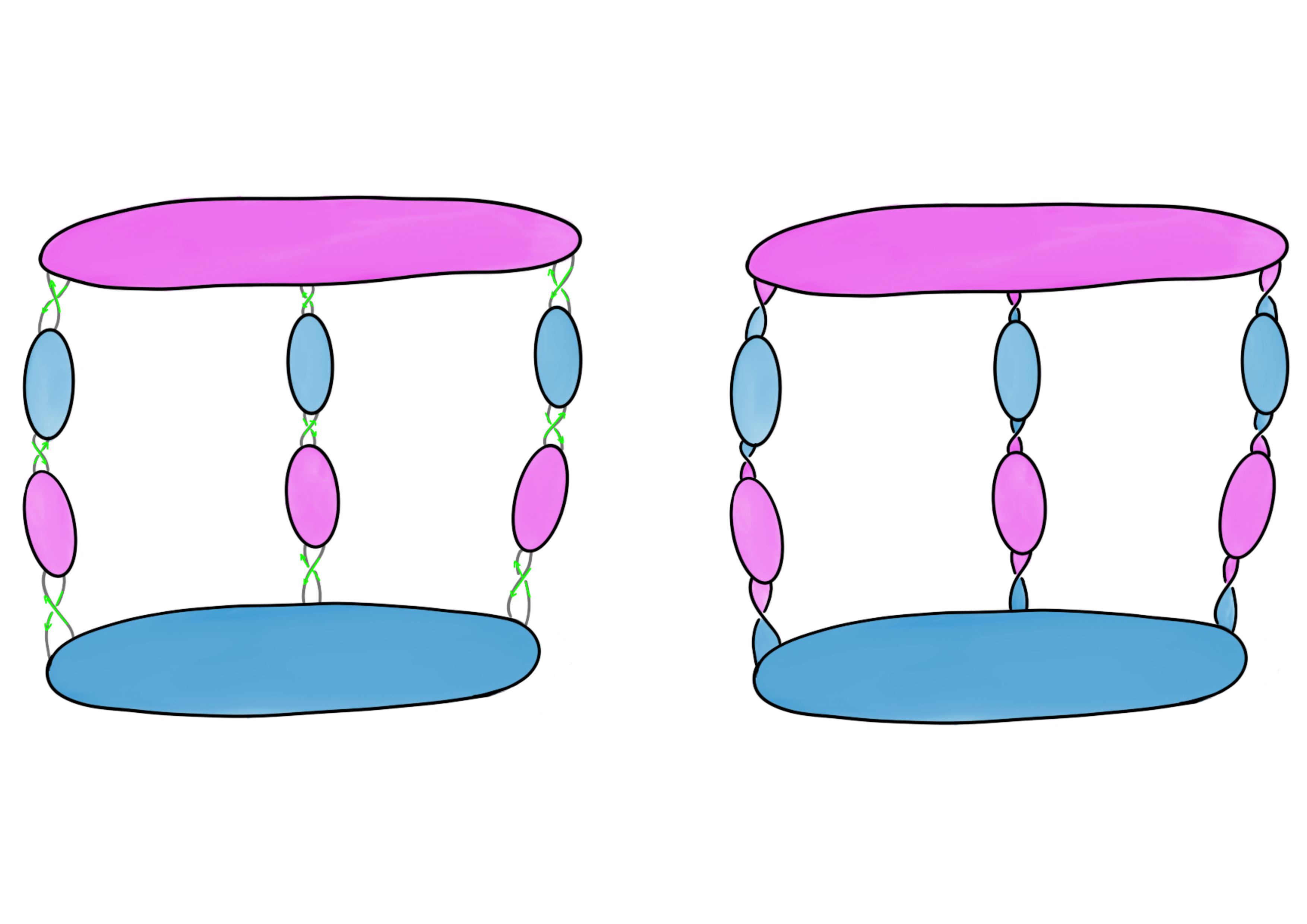}
    \caption{Superficie final}
    \label{seifertsup3}
\end{figure}

Note que en la \Cref{seifertsup3} usamos el color rosado para los discos cuya frontera est\'a orientada en el sentido de las manecillas del reloj, y azul auqellos con la orientaci\'on contraria. Por \'ultimo para reincorporar los cruces del nudo a nuestra superficie unimos cada disco mediante bandas trenzadas, el giro que debe tener la banda se infiere de la orientaci\'on original que ten\'ia el cruce del nudo.\\

\begin{figure}[H]
    \centering
    \includegraphics[clip, trim=0 0 0 0, width=0.85 \linewidth]{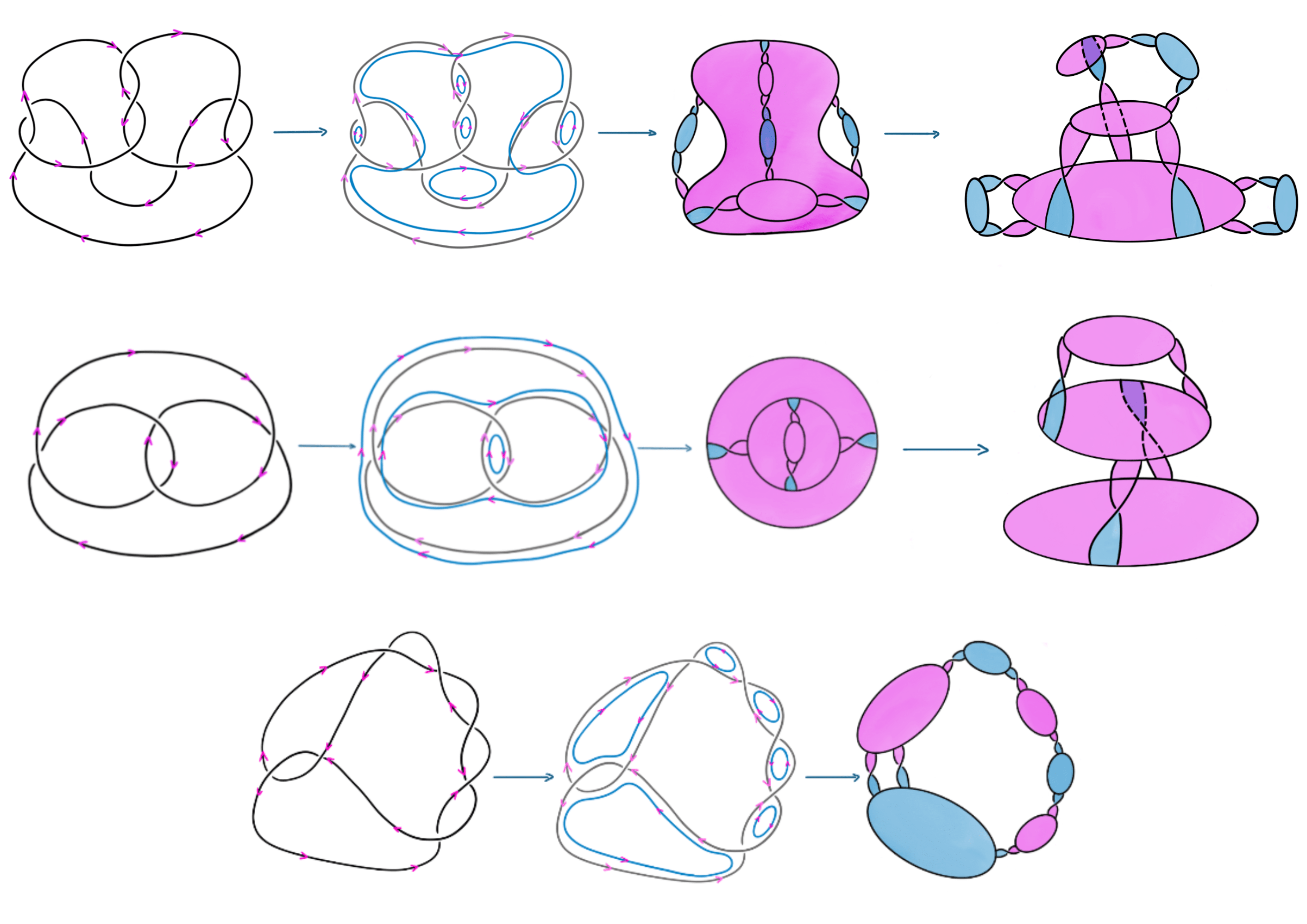}
    \caption{Superficies de Seifert de distintos nudos}
    \label{seifertejemplos}
\end{figure}

En el ejemplo anterior no se considera un nudo en donde tengamos c\'irculos anidados. Para ilustrar la separaci\'on de discos anidados usando la altura de anidaci\'on, inclu\'imos la \Cref{seifertejemplos}. 

\begin{figure}[H]
    \centering
    \includegraphics[clip, trim=30 800 900 100, width=0.4\linewidth]{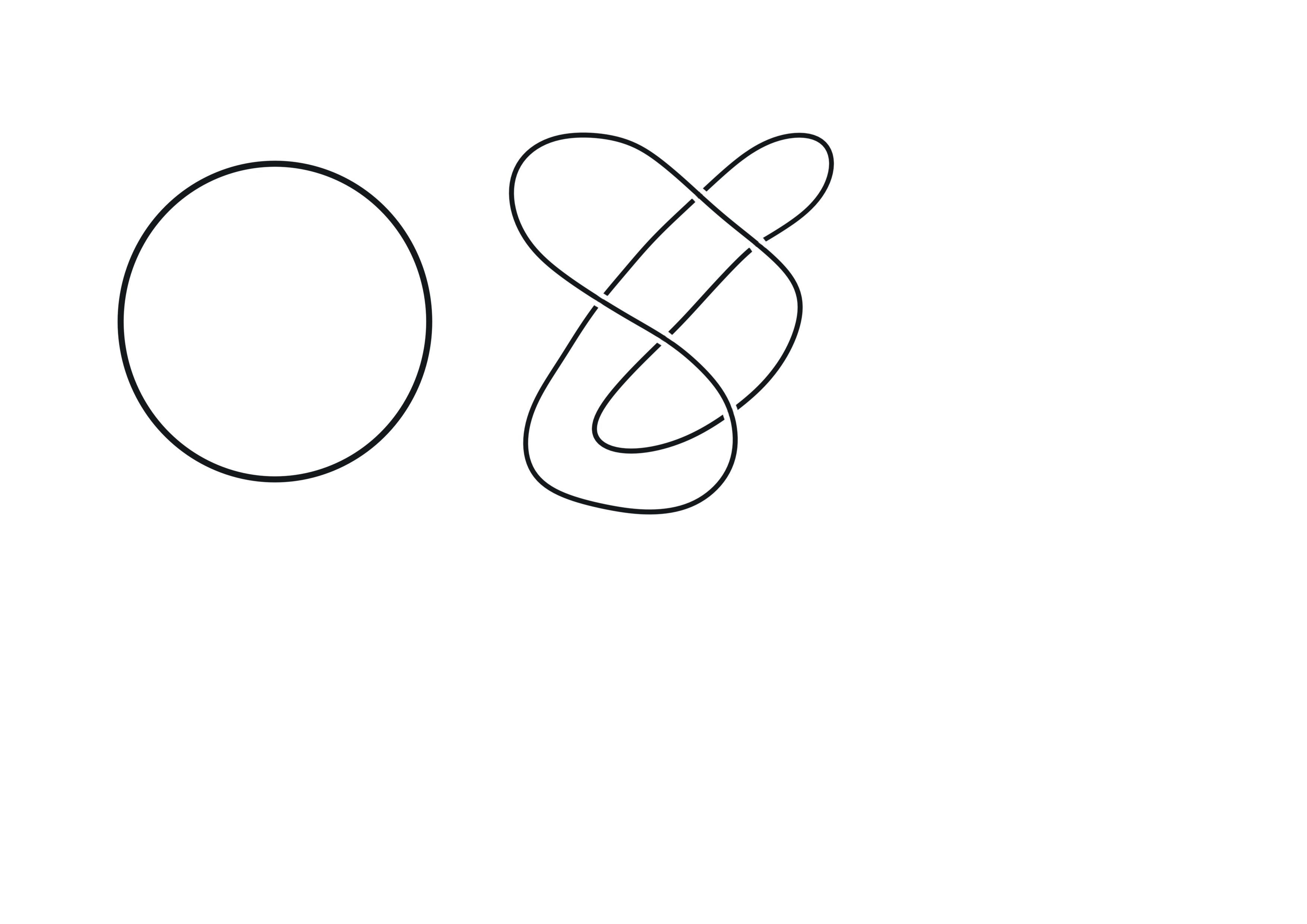}\includegraphics[clip, trim=30 300 30 150, width=0.4\linewidth]{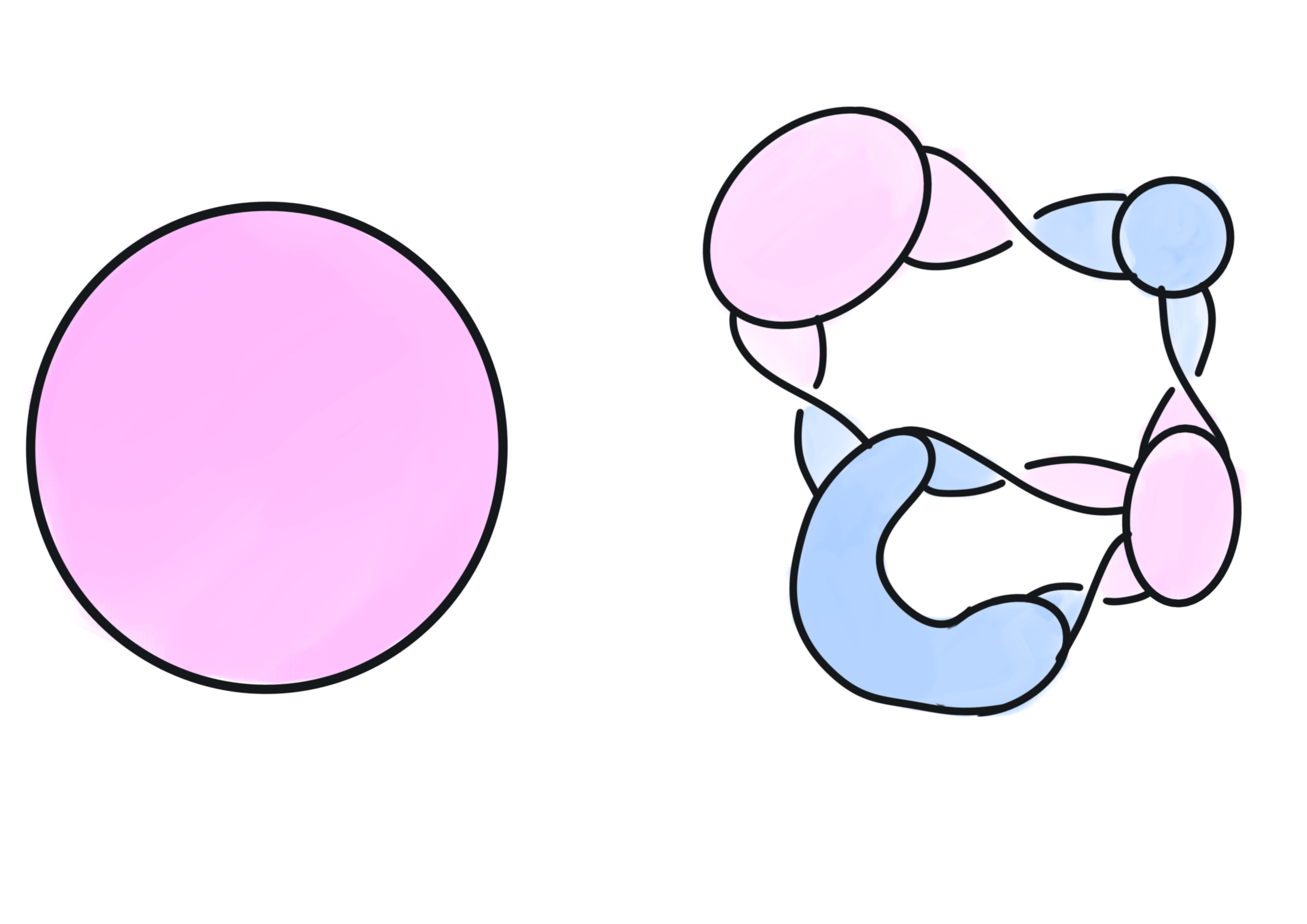}
    \caption{Dos diagramas diferentes del nudo trivial y las superficies obtenidas a partir del algoritmo de Seifert. Sus caracter\'isticas de Euler son $1$ y $-1$ respectivamente.}
    \label{diag-trivial}
\end{figure}

Aunque este algoritmo nos permite construir la superficie de Seifert asociada, su resultado depende del diagrama del nudo, y por tanto no siempre nos da como output la superficie de g\'enero m\'inimo. Por ejemplo, la \Cref{diag-trivial} muestra dos diagramas del nudo trivial y las superficies obtenidas al `correr' el algoritmo de Seifert usando cada diagrama como `input'. De acuerdo a el \Cref{clas-frontera}, para distinguir las dos superficies ser\'ia suficiente calcular su caracter\'istica de Euler. Nuestro problema en \'este momento es que las f\'ormulas que tenemos para calcularla est\'an dadas en t\'erminos de v\'ertices, aristas, y caras \cref{euler} o en t\'erminos del g\'enero \cref{euler-g-frontera}. La siguiente proposici\'on nos da una f\'ormula que podremos aplicar directamente a las superficies obtenidas como `output' del algoritmo de Seifert. 

\proposicion{}{Si una superficie $F\subset\mathbb{R}^3$ compacta, orientable y que en su frontera tiene una componente conexa, se construye a partir de $b$ bandas y $d$ discos, entonces 
\[
\chi(F)=d-b
\]} 

\demostracion{
El resultado del algoritmo de Seifert se puede pensar como el espacio cociente de la uni\'on de los discos de Seifert y las bandas, con relaci\'on de equivalencia determinada por la identificaci\'on de los lados `cortos' de las bandas con arcos en los c\'irculos de Seifert. 

La uni\'on consiste de $b+d$ caras, $4b$ aristas y $4b$ v\'ertices de las bandas, y adicionalmente $2b$ aristas y $2b$ v\'ertices de los discos. Para ver \'esto \'ultimo, note que cada c\'irculo de Seifert $C$ se construye como la uni\'on de un conjunto de arcos provenientes del nudo $K$, y otro proveniente de las resoluciones de los cruces. Puesto que $\chi(C)=0$, \'estos dos conjuntos tienen el mismo n\'umero de elementos que se denota por $a_C$. La suma de $a_C$ sobre todos los c\'irculos de Seifert $C$ es igual al n\'umero total de cruces del diagrama, y \'este a su vez es igual al n\'umero $b$ de bandas. 

En el espacio cociente el n\'umero de caras sigue siendo $C=b+d$. Sin embargo, dado que los lados `cortos' de las bandas se identifican con arcos en los c\'irculos de Seifert, el n\'umero de aristas y v\'ertices es menor. De hecho, todos los v\'ertices de los c\'irculos de Seifert se identifican con v\'ertices en las bandas por lo que $V=2b$. A su vez, la mitad de las aristas de las bandas se identifican con aristas de los c\'irculos de Seifert y por lo tanto $A=2b+2b$. Como consecuencia, $$\chi=(2b)-(2b+2b)+(b+d)=d-b.$$
}

\section{Invariantes de Nudos}\label{invariantes}
A grandes rasgos, un invariante es una aplicaci\'on que asigna a cada uno un objeto matem\'atico familiar -y usualmente algebr\'aico. Debe asignar el mismo objeto a nudos equivalentes. En \'esta secci\'on nos concentraremos en invariantes de nudos que se definen en t\'erminos de las superficies de Seifert.

\subsection{El 3-g\'enero de un nudo}
Como hemos visto, las superficies se clasifican en t\'erminos de su g\'enero (\Cref{clasificacion}, \Cref{clas-frontera}), y que a cada nudo $K$ podemos asociar una superficie $F$ con el nudo $K$ como frontera~(\Cref{seifert}). Es entonces tentador asignar la cantidad $g(F)$ al nudo $K$ y as\'i definir un invariante de nudos. El problema es que esta asignaci\'on est\'a lejos de estar bien definida como lo ilustra \Cref{diag-trivial} donde se muestran superficies de g\'enero $0$ y $1$ con frontera el nudo trivial. La soluci\'on es simple: tomar el m\'inimo sobre todas las posibles superficies con frontera el nudo dado. 

\dfe{3-g\'enero de un nudo}{ 
Sea $K\subset\mathbb{R}^3$, se define el g\'enero del nudo $K$ como
$$g_3(K)=\min\{g(F) \mid F\subset\R^3 \text{ es una superficie } \partial F=K\}.$$
}

Dado que se consideran \textit{todas} las superficies con frontera un nudo dado $K$, de la misma definici\'on se sigue que el $3$-g\'enero es un invariante de nudos, es decir, si dos nudos son isot\'opicos, entonces sus $3$-g\'enero son iguales. Lo que no es cierto es que si dos nudos tienen el mismo $3$-g\'enero, entonces los nudos son equivalentes. \'Esto lo veremos m\'as adelante cuando consider\'emos el nudo tr\'ebol y el nudo en ocho. Sin embargo, el caso del nudo trivial es especial: 

\teorema{3-G\'enero del nudo trivial}{\label{g3-trivial}
Un nudo $K\subset\mathbb{R}^3$  satisface $g_3(K)=0$ si y solamente si $K$ es el nudo trivial.}

La importancia de este teorema radica en poder reconocer directamente con el g\'enero si un nudo es isot\'opico al nudo trivial, razonamiento que no se puede establecer con los dem\'as nudos. \\

\subsection{Matriz de Seifert} 
Para obtener invariantes de nudos a partir de superficies de Seifert, vamos a extender la clasificaci\'on de las superficies con frontera conexa y nos vamos a concentrar no s\'olo en el tipo de homeomorfismo de una superficie $F\subset\R^3$, sino tambi\'en en la manera precisa en la que $F$ est\'a encajada. Para ello cuantificaremos la interacci\'on entre curvas cerradas en $\R^3\setminus F$ y en $F$. Por ejemplo, si tomamos una superficie de Seifert para un nudo construida a partir de discos y bandas, podemos tomar una familia de curvas orientadas como en la \Cref{seifert+base} que pasen por el centro de las bandas y formen arcos dentro de los discos. La manera en que estas curvas se retuercen y enlazan lleva informaci\'on sobre el nudo, la cual se captura dentro de la \emph{matriz de Seifert del nudo}, a partir de la noci\'on del  \emph{n\'umero de enlace} y el \emph{push-off} de curvas en $F$. \\

\begin{figure}[H]
    \centering
    \includegraphics[width=0.3\linewidth]{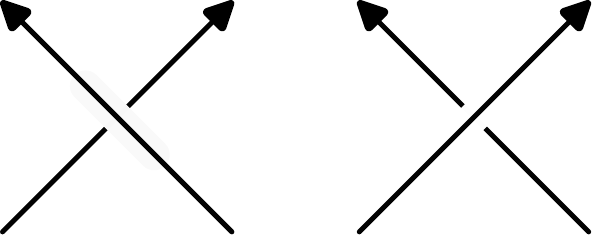}
    \caption{Cruce de la mano-izquierda ($-1$) y de la mano-derecha ($+1$).}
    \label{cruces}
\end{figure}

Sean $\alpha$ y $\beta$ dos curvas cerradas y orientadas en $\R^3$, el \emph{n\'umero de enlace} de  $\alpha$ y $\beta$ se denota por $lk(\alpha,\beta)$ y se calcula de la siguiente manera: 
\begin{enumerate}
\item Sea $D$ un diagrama de la uni\'on $\alpha\cup\beta$.
\item Se miran los cruces en $D$ y se consideran solamente aquellos en los que una hebra de $\alpha$ pasa por encima de una hebra de $\beta$.
\item A cada uno de los cruces se le asigna un $-1$ si modulo rotaciones el cruce se ve como el cruce a la izquierda de la \Cref{cruces}. De lo contrario se asigna un $+1$.
\item El n\'umero de enlace es la suma de los n\'umeros obtenidos en el paso anterior.
\end{enumerate}

\nota{}{ Propiedades del n\'umero de enlace.  
\begin{itemize}
    \item El n\'umero de enlace no depende del diagrama, sino solo de las orientaciones de las componentes.
    \item El n\'umero de enlace siempre es un entero, no cambia bajo deformaciones y es sim\'etrico: \( lk(\alpha, \beta) = lk(\beta, \alpha) \).
    \item La suma se puede separar en dos partes: los cruces donde \( \alpha \) est\'a por encima de \( \beta \) y los cruces donde \( \beta \) est\'a por encima de \( \alpha \).
\end{itemize}}

\begin{figure}[H]
    \centering
    \includegraphics[width=0.7\linewidth]{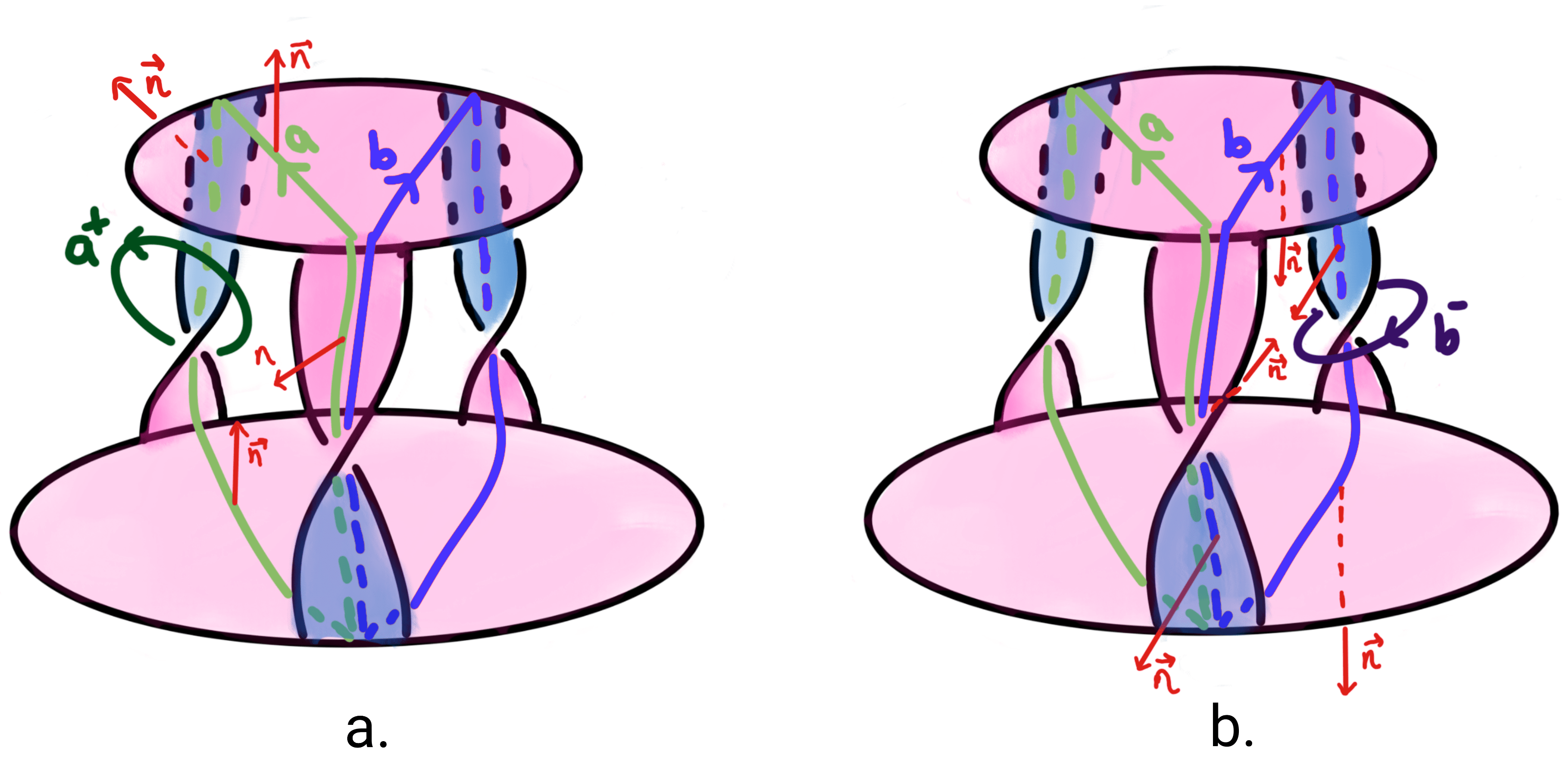}
    \caption{Push off positivo y negativo de dos curvas}
    \label{push offs}
\end{figure}

Para definir el push-off de $\gamma$ una curva en $F\subset \R^3$ una superficie, primero identificamos una vecindad cerrada de $F$ con el producto $F\times [-1,1]$. El push-off positivo $\gamma^+$ de $\gamma$ es la copia de $\gamma$ en $F\times\{1\}$ de modo que $\gamma^+$ se encuentre justo por `encima' de $F$. Como la superficie $F$ es orientable, es posible escoger un campo vectorial unitario $\vec{n}$ que es normal a $F$ y que var\'ia continuamente en cada punto. Desde \'esta perspectiva $\gamma^+$ es la curva en $\R^3\setminus F$ que corre en paralelo a $\gamma$ y se obtiene al \textit{``halar''} a $\gamma$ el la direcci\'on de $\vec{n}$. Para obtener el push-off negativo $\gamma^-$ usamos la normal negativa $-\vec{n}$ (ver Figura \ref{push offs}).

\dfe{Matriz de Seifert}{
Sea $F$ una superficie de Seifert para un nudo $K$, y $a_1,a_2,\ldots,a_{2g}$ curvas cerradas en $F$ que representan una base para $H_1(F;\Z)\cong\Z^{2g}$. La matriz $V$ de tama\~no $2g\times 2g$ con entradas $$v_{ij}=lk(a_i,a_j^+)$$ es una matriz de Seifert para $K$ }

Consideremos las siguientes superficies de Seifert, junto con sus curvas y sus respectivos push off positivos:

\begin{figure}[H]
    \centering
    \includegraphics[width=\textwidth]{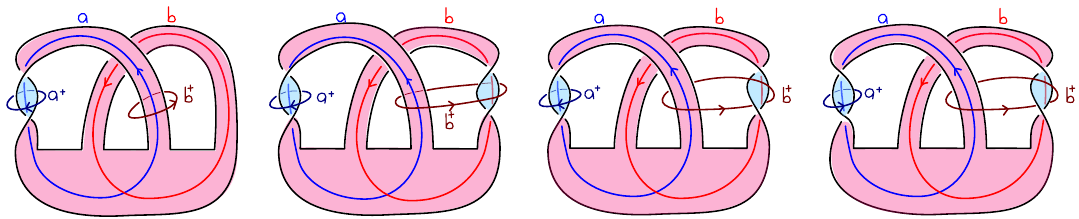}\\
    \text{Entonces, sus matrices de Seifert asociadas respectivamente son:}\\
    
    \begin{center}
$V_1=\left[ \begin{array}{rr}
   	1 & 1\\
        	0 & 0
    \end{array}\right]$ \hfill $V_2= \left[ \begin{array}{rr}
	1 & 1\\
	0 & 1
    \end{array}\right]$\hfill $V_3= \left[ \begin{array}{rr}
	1 & 1\\
	0 & -1
    \end{array}\right]$\hfill $V_4= \left[ \begin{array}{rr}
-1 & 1\\
 0 & -1
\end{array}\right]$
\end{center}

    \caption{Ejemplos de superficies de Seifert, base de $H_1(F;\Z)$ y sus respectivos `push-offs'.}
    \label{seifert+base}
\end{figure}

Una matriz de Seifert depende primero de la elecci\'on de superficie $F$ y segundo de la elecci\'on de una base particular de $H_1(F;\Z)$, por lo que, por s\'i sola, no es un invariante. Para obtener invariantes de nudos debemos primero establecer la relaci\'on precisa entre matrices de Seifert del mismo nudo, pero obtenidas a partir de elecciones distintas superficies o de bases de homolog\'ia. Dos matrices de Seifert $V_1,V_2$ asociadas a una misma superficie $F$ pero a bases distintas de $H_1(F;\Z)$ son congruentes de forma unimodular. Es decir, existe una matriz unimodular $P$ tal que $V_2=P^{T}V_1P$. La relaci\'on entre matrices de Seifert asociadas a superficies diferentes es un poco m\'as delicada y depende de la siguiente definici\'on. 

\begin{definition} Sea $V\in M_{n\times n}(\Z)$. Una matriz $W$ es una ampliaci\'on elemental de $V$ si $W$ tiene la forma \begin{center}
\begin{tabular}{ccc}
$W=\left[\begin{array}{c|cc}V & 0 & 0 \\\hline a & 0 & 0 \\ 0 & 1 & 0\end{array}\right]$ & o &
$W=\left[\begin{array}{c|cc}V & b & 0 \\\hline 0 & 0 & 1 \\ 0 & 0 & 0\end{array}\right]$
\end{tabular}\end{center}
para algunos vectores  $a\in M_{1\times n}(\Z)$ y $b\in M_{n\times 1}(\Z)$. En \'este caso $V$ se llama una reducci\'on elemental de $W$.
\end{definition}

\dfe{$S$-equivalencia}{\label{S-equiv}
Dos matrices $V_1$ y $V_2$ son $S$-equivalentes si una puede obtenerse de la otra mediante un n\'umero finito de ampliaciones/reducciones elementales, y congruencia unimodular.
}

\teorema{}{\label{S-inv}
Sea $K$ un nudo. 
Si $V_1$ y $V_2$ son dos matrices de Seifert asociadas al mismo nudo $K$, entonces 
$V_1$ y $V_2$ son $S$-equivalentes. Es m\'as, en \'este caso tenemos 
\begin{enumerate}[label=(\arabic*)]
\item\label{det} Existe $m\in\Z$ tal que $\det{\left(V_2-tV_2^T\right)}=(\pm t^m)\cdot\det{\left(V_1-tV_1^T\right)}$, para todo $t\in\R$.
\item\label{sig} $\sigma\left(V_2+V_2^T\right)=\sigma\left(V_1+V_1^T\right)$
\end{enumerate}
}

\demostracion{A grandes rasgos, la demostraci\'on se basa en las siguientes ideas:
\begin{enumerate}[label=(\alph*)]
\item Las ampliaciones elementales provienen de realizar $0$-cirug\'ia a superficies,
\item las reducciones elementales provienen de realizar $1$-cirug\'ia a superficies, y
\item la congruencia unimodular representa cambios de base.
\end{enumerate}

La demostraci\'on completa requiere un estudio cuidadoso del efecto de la cirug\'ia en $H_1(F;\Z)$ y en la forma de Seifert. Puede encontrarse en \cite[Teorema 8.4]{lickorish}.\\

Para establecer el resultado es suficiente considerar la congruencia unimodular y los dos tipos de ampliaciones de manera independiente. Primero, dos matrices $V_1$ y $V_2$ se relacionan por congruencia unimodular si existe una matriz unimodular $P$ con entradas en $\Z$ tal que $V_2=PV_1P^T$. En ese caso
$$\det{\left(V_2-tV_2^T\right)}=\det{(P\left(V_1-tV_1^T\right)P^T)}=\det{(P)}^2\det{\left(V_1-tV_1^T\right)}=\det{\left(V_1-tV_1^T\right)}, $$ y $\sigma\left(P\left(V_2+V_2^T\right)P^T\right)=\sigma\left(V_1+V_1^T\right)$.\\

Para ver que el resultado es válido para los casos restantes, sin pérdida de generalidad supongamos $$\fontsize{9}{7}\selectfont V_2=\left[\begin{array}{c|cc}V_1 & 0 & 0 \\\hline a & 0 & 0 \\ 0 & 1 & 0\end{array}\right] $$
entonces tenemos
$$\fontsize{9}{7}\selectfont V_2-tV_2^T=\left[\begin{array}{c|cc}V_1-tV_1^T & -ta^T & 0 \\\hline a & 0 & -t \\ 0 & 1 & 0\end{array}\right], $$
y usando la última fila para calcular el determinante vemos que $$\fontsize{9}{7}\selectfont \det{\left(V_2-tV_2^T\right)}=\det{\left[\begin{array}{c|c}V_1-tV_1^T  & 0 \\\hline a & -t \end{array}\right] }=(-t)\cdot\det{\left(V_1-tV_1^T\right)}.$$

Adicionalmente, 
$$\fontsize{9}{7}\selectfont 
V_2+V_2^T=\left[\begin{array}{c|cc}V_1+V_1^T & a^T & 0 \\[1ex]\hline a & 0 & 1 \\ 0 & 1 & 0\end{array}\right].$$
 Sea $I$ la matriz identidad con las mismas dimensiones de $V_1$. Si $Q$ es la matriz 
$$\fontsize{9}{7}\selectfont 
Q=\left[\begin{array}{c|cc}I_n & 0 & -a^T \\[2ex]\hline 0 & 1 & 0 \\[1ex] 0 & 0 & 1\end{array}\right],$$
entonces
$$\fontsize{9}{7}\selectfont Q\left(V_2+V_2^T\right)Q^T=\left[\begin{array}{c|cc}V_1+V_1^T & 0 & 0 \\\hline 0 & 0 & 1 \\ 0 & 1 & 0\end{array}\right].$$
Dado que $\fontsize{9}{7}\selectfont \sigma\left(\left[\begin{array}{cc}0 & 1 \\ 1 & 0\end{array}\right]\right)=0$ y la signatura es aditiva, tenemos $$\sigma\left(Q\left(V_2+V_2^T\right)Q^T\right)=\sigma\left(V_1+V_1^T\right).$$ 
}

\subsection{Polinomio de Alexander y Signatura}
El polinomio de Alexander fue uno de los primeros invariantes de nudos en ser definidos. La definici\'on original de Alexander se obtiene desde la combinatoria y data de 1923. En los 100 a\~nos desde su introducci\'on, decenas de procedimientos para calcular el polinomio se han desarrollado, inclu\'idas las relaciones de Skein \cite[pp. 207-215]{livingston}, \'el c\'alculo de Fox \cite{fox} y la homolog\'ia de espacios recubridores \cite{milnor}. En estas notas nos enfocaremos en la definici\'on algebr\'aica a partir de la matriz de Seifert. Para definir la signatura de un nudo, record\'emos que dada una matriz sim\'etrica, su signatura es la resta del n\'umero de valores propios positivos menos el n\'umero de valores propios negativos, y que para cualquier matriz $V$, la matriz $V+V^T$ es sim\'etrica.

\begin{definition}\label{invariants} Sea $K\subset\R^3$ con matriz de Seifert $V$. 
\begin{enumerate}[label=(\arabic*)]
\item El polinomio de Alexander de $K$ se define como $$\Delta_K(t)=\det{\left(V-tV^T\right)}$$ (m\'odulo m\'ultiplos de $\pm t^m$)
\item La signatura de $K$ es la signatura de la matriz sim\'etrica $V+V^T$. Se denota como $\sigma(K)$
\end{enumerate}
\end{definition}

Para demostrar que el polinomio de Alexander y la signatura son efectivamente invariantes de nudos, se debe demostrar que sus valores no dependen de la elecci\'on de superficie y/o matriz de Seifert. El siguiente teorema es simplemente una reformulación de \Cref{S-inv}.

\teorema{}{
Si $K_1$, $K_2$ son nudos equivalentes,
\begin{enumerate}[label=(\arabic*)]
\item (\cite[Ch. 6 Sec. 2 Cor. 4]{livingston}) $\Delta_{K_1}(t)=\Delta_{K_2}(t)$ (m\'odulo m\'ultiplos de $\pm t^m$)
\item (\cite[Ch. 6 Sec. 3 Cor. 5]{livingston}) $\sigma(K_1)=\sigma(K_2)$
\end{enumerate}
}

Detr\'as de \'estas definiciones y teoremas se esconde una pregunta: si el polinomio de Alexander y la signatura se derivan de una superficie de Seifert, pero son independientes de la elecci\'on espec\'ifica de \'esta, ¿se pueden relacionar con el 3-g\'enero 3 del nudo? La respuesta es un rotundo s\'i: el polinomio de Alexander y la signatura permiten establecer cotas inferiores para el 3-g\'enero. Para ello, definimos $\deg(\Delta_K(t))$ como la diferencia entre el mayor y el menor de los grados de los monomios del polinomio. N\'otese que \'esta cantidad no cambia incluso m\'odulo m\'ultiplos de de potencias de la variable $t$.

\teorema{Cotas inferiores al 3-g\'enero}{\label{cotas}
Sea $K\subseteq S^3$ un nudo. Se tienen las siguientes cotas para el 3-g\'enero de $K$
\begin{enumerate}
	\item\label{delta-g} $\deg(\Delta_K(t))\leq 2g_3(K)$
	\item\label{sigma-g} $|\sigma(K)|\leq 2g_3(K)$
\end{enumerate}
}

La demostraci\'on se obtiene casi directamente de las siguientes ideas que se estudian en un primer curso de algebra lineal: Sea $X$ una matriz cuadrada con entradas $x_{i,j}$ $1\leq i,j\leq n$. 
\begin{enumerate}
    \item La f\'ormula de Leibniz expresa $\det(X)$ como la suma $$\sum _{\tau \in S_{n}}\operatorname {sgn}(\tau )\prod _{i=1}^{n}x_{\tau (i),i}\footnote{En \'esta expresi\'on $S_n$ denota el grupo de permutaciones de un conjunto de $n$ elementos, y $\operatorname {sgn}(\tau )$ representa el signo de $\tau$}.$$ Esta expresi\'on nos permite considerar $\det(X)$ como un elemento de $\Z[x_{i,j}]$ de grado a lo sumo $n$.
    \item Supongamos que $X$ es sim\'etrica, y $n_+,n_-$ denotan respectivamente el n\'umero de valores propios positivos y negativos de $X$. Entonces $\sigma(X)=n_+-n_-\leq n_++n_-\leq n$.
\end{enumerate}

\demostracion{
Sea $F$ una superficie en $\R^3$ de g\'enero $g$ y con frontera el nudo $K$. Sea $V$ una matriz de Seifert asociada a $F$. Recordemos que $V$ es entonces una matriz de dimensi\'on $2g\times 2g$. Dado que el determinante es una combinaci\'on lineal de productos de $2g$ entradas de la matriz, el grado de $\det(V-tV^T)$ es a lo sumo $2g$. De la misma manera, $\sigma(K)\leq 2g$. Escogiendo entonces la superficie $F$ de tal manera que $g(F)=g_3(K)$ obtenemos las cotas deseadas. }

\subsection{Ejemplos}
   Observemos algunos ejemplos del c\'alculo de polinomios de Alexander. Para ello consideremos las superficies y matrices de Seifert de \Cref{seifert+base}. Con estas matrices consideramos las matrices $V-tV^T$, y sus respectivos determinantes que dan lugar a los polinomios de Alexander correspondientes a cada superficie, como se muestran a continuaci\'on.

\begin{table}[h]
\renewcommand{\arraystretch}{2}
\newcolumntype{C}{>{\centering\arraybackslash $}X<{$}}%
\begin{tabularx}{0.75\textwidth}{|>{$}l<{$}| C |C|C|C|}
\hline
K& U&T_{2,-3}& 4_1& T_{2,3}\\
\hline%
\Delta_K(t)& t & 1-t+t^2 & -1+3t-t^2 & 1-t+t^2 \\
\hline
\sigma(K)& 0& 2 & 0 & -2\\
\hline
\end{tabularx}
\caption{Ejemplos de c\'alculos del polinomio de Alexander y la signatura de los nudos de la \Cref{seifert+base}}
\label{ejemplos-invariantes}
\end{table}


Si bien el 3-g\'enero $g_3(K)$ de un nudo $K$ es un invariante, su definici\'on como el m\'inimo de un conjunto dificulta su c\'alculo. En algunos casos el \Cref{cotas} nos da herramientas suficientes para calcular valores exactos de $g_3(K)$, como lo describiremos a continuaci\'on.\\

Consider\'emos el nudo tr\'ebol $T$. En \Cref{ejemplos-invariantes} se muestra que $|\sigma(T)|=2$ y $deg(\Delta_T(t))=2$. Como consecuencia de \Cref{cotas} tenemos entonces $1\leq g_3(T)$. Por otro lado, \Cref{seifert+base} muestra una superficie de Seifert de $T$ con g\'enero exactamente 1, de lo cual se sigue que $g_3(T)=1$. \\

El caso del nudo en ocho $4_1$ es similar: $|\sigma(4_1)|=0$, $deg(\Delta_{4_1}(t))=2$ (ver \Cref{ejemplos-invariantes}) y  \Cref{seifert+base} muestra una superficie de Seifert para $4_1$ con g\'enero exactamente 1. Si bien la signatura nos da una cota inferior para $g_3(4_1)$, no nos permite decidir si $g_3(4_1)$ es $0$ o $1$. Sin embargo, el grado del polinomio de Alexander nos permite concluir que $g_3(4_1)=1$. \'Este ejemplo ilustra la importancia de calcular ambas cotas. \\


\section{Concordancia}
En esta secci\'on se quiere abordar la estructura de grupo que adoptan los nudos a partir de una relaci\'on de equivalencia, siendo la suma conexa la operaci\'on del grupo. El objetivo de la secci\'on es dar una breve introducci\'on a esta estructura y motivar a quien lea estas notas para indagar m\'as sobre este tema.\\ 

\subsection{Monomio de grupo}
\begin{figure}[H]
    \centering
    \includegraphics[width=0.5\linewidth]{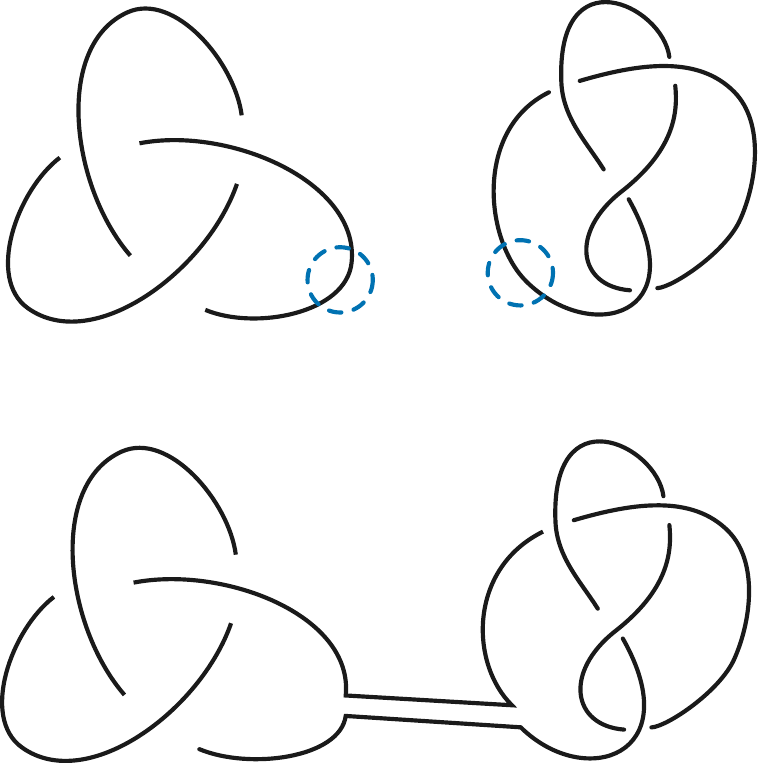}
    \caption{Suma conexa}
    \label{suma conexa}
\end{figure}

Sean $K_1$ y $K_2$ dos nudos. Para entender de manera intuitiva la suma conexa $K_1\# K_2$, primero se toman proyecciones disjuntas de $K_1$ y $K_2$ en el plano. Luego se escoge un arco $a$ con frontera $\partial a=\{a_1,a_2\}$ tal que $a_i\in K_i$ y $\text{int}(a)$ es disjunto de los nudos. Si $U_i$ denotan vecindades de $a_i$ en $K_i$, la suma conexa es el resultado de unir $K_1\setminus U_1$ y $K_2\setminus U_2$ a lo largo de arcos paralelos a $a$, como se muestran en la \Cref{suma conexa}. La suma conexa hereda una orientaci\'on consistente con las orientaciones de los dos nudos originales, y la clase de isotop\'ia es independiente de la elecci\'on de arco $a$. 

Para estudiar la manera en que nuestros invariantes se comportan bajo la suma conexa, es conveniente introducir la noci\'on de \emph{suma conexa en la frontera} de superficies de Seifert. Para ello, si $F_i$ es una superficie de Seifert para $K_i$ ($i=1,2$), se construye la superficie $F_1\natural F_2$ como el espacio cociente de la uni\'on de $F_1$, $F_2$ y una vecindad $a\times [0,1]$ del arco $a$ que identifica los segmentos `verticales' $a_i\times [0,1]$ con las vecindades $U_i\subset F_i$. Como consecuencia de la construcci\'on misma tenemos que la frontera de $F_1\natural F_2$ es el nudo $K_1\# K_2$, y adicionalmente
\begin{equation}
g_3(K_1\# K_2)\leq g(F_1\natural F_2)=g(F_1)+g(F_2). 
\end{equation}
Tomando m\'inimos sobre todas las posibles superficies $F_1$ y $F_2$ obtenemos la desigualdad $g_3(K_1\#K_2)\leq g_3(K_1)+g_3(K_2)$. La demostraci\'on de la desigualdad contraria se puede encontrar en \cite[5.A.14]{rolfsen}.\\

La siguiente proposici\'on establece el comportamineto preciso de los invariantes estudiados en la \Cref{invariantes} con respecto a la suma conexa.

\proposicion{}{Sean $K_1,K_2$ nudos. Entonces,
\begin{enumerate}
\item $g_3(K_1\#K_2)=g_3(K_1)+g_3(K_2)$
    \item $\Delta_{K_1\#K_2}(t)=\Delta_{K_1}(t)\cdot\Delta_{K_2}(t)$.
    \item $\sigma(K_1\# K_2)=\sigma(K_1)+\sigma(K_2)$.
\end{enumerate}}

Los \'ultimos dos enunciados de esta proposici\'on se demuestran eligiendo primero conjuntos $B_1$ y $B_2$ abiertos disjuntos de $\R^3$ que contienen a los nudos $K_1$ y $K_2$ respectivamente, y luego superficies de Seifert $F_1$ y $F_2$ completamente contenidas en $B_1$ y $B_2$ respectivamente. La matriz de Seifert asociada a $F_1\natural F_2$ en este caso es una matriz por bloques de la forma 
$$
\begin{bmatrix}
   V_1 & 0\\
    0 & V_2
\end{bmatrix},
$$
donde $V_i$ es la matriz de Seifert asociada a $F_i$ ($i=1,2$). \\

Hasta el momento hemos establecido que el conjunto de nudos orientados forma un semigrupo conmutativo con suma conexa como operaci\'on binaria. Para obtener una estructura de grupo tendr\'iamos que demostrar la existencia de un elemento identidad, y la existencia de inversos. No es dif\'icil ver que el nudo trivial es el elemento identidad, sin embargo, como el siguiente lema lo demuestra, los nudos no triviales no tienen inversos.

\lema{}{ El nudo trivial no se puede expresar como la suma de dos nudos no triviales.}
\demostracion{ Si $K$ y $J$ son dos nudos no triviales tales que $K\#J$ es el nudo trivial $U$, entonces $g_3(K\#J)=0$, y $g_3(K),\,g_3(J)>0$ como consecuencia de \Cref{g3-trivial}. Dado que el 3-g\'enero es aditivo, tenemos que $0=g_3(K\#J)=g_3(K)+g_3(J)>0$,  llegando así a una contradicción.}

\subsection{Grupo de Concordancia}
La falta de inversos se puede resolver considerando superficies encajadas ya no en $\R^3$, sino en $\R^4$, y con ellas la noci\'on de concordancia de nudos. El estudio de la concordancia de nudos comenz\'o con el trabajo de Artin~\cite{artin} quien not\'o que dado un encaje topol\'ogico diferenciable $\phi: S^2\hookrightarrow \R^4$, la intersecci\'on de su imagen con el subespacio $\R^3\times\{0\}$ es un nudo $K$ no trivial, y que $K$ es la frontera de un disco propiamente encajado en $\R^4$. Los nudos que surgen de esta manera se conocen como nudos \emph{`rebanada' o `tajada'}. Para cuantificar \'esta idea, extendemos la definici\'on de $g_3$ de la siguiente manera:

\dfe{4-g\'enero de un nudo}{\label{g4} 
Sea $K\subset\mathbb{R}^3$, se define el g\'enero del nudo $K$ como
$$g_4(K)=\min\{g(F) \mid \phi: F\hookrightarrow \R^4_+ \text{ encaje diferenciable propio } \phi\left(\partial F\right)=K\}.$$}

En la definici\'on anterior $\R^4_+$ denota el conjunto de puntos de $\R^4$ cuya \'ultima coordenada es mayor o igual a cero, de tal manera que $\R^3\times\{0\}$ es su frontera. 

Con miras a una definici\'on rigurosa de concordancia, definimos:

\dfe{Imagen Especular}{Sea $m: \R^3\to \R^3$ la funci\'on definida por $(x,y,z)\to (x,y,-z)$. Dado un nudo $K$, su imagen especular es el nudo obtenido como $m(K)$}

En t\'erminos de diagramas de nudos, si $D$ es el diagrama de $K$ que se obtiene a partir de la proyecci\'on de $K$ en el plano $xy$, entonces el diagrama que se obtiene al cambiar en cada cruce de $D$ las hebras superiores por hebras inferiores es un diagrama de $m(K)$. N\'otese que si $K$ es un nudo orientado, $m(K)$ es tambi\'en un nudo orientado con orientaci\'on inducida por aquella de $K$.\\

Dado un nudo orientado $K\subset S^{3}$, denotamos por $-K$ el nudo $m(K)$ con la orientaci\'on opuesta a la de $K$. 

\dfe{Concordancia}{\label{conc}
Dos nudos $K$ y $J$ en $S^3$ son concordantes si $g_4(K\#-J)=0$. En otras palabras, si la suma conexa $K\#-J$ es la frontera de un disco $D^2$ propiamente y suavemente encajado en $\R^4_+$.}

Vemos entonces la manera en que la noci\'on de nudo tajada da lugar a la idea de concordancia, y con ella logramos establecer la estructura de grupo que busc\'abamos:

\teorema{Fox-Milnor~\cite{fox-milnor}}{ La concordancia es una relaci\'on de equivalencia y el conjunto de clases de concordancia forma un grupo abeliano $\mathcal{C}$ con la suma conexa como operaci\'on binaria.}

El elemento trivial de $\mathcal{C}$ es la clase de equivalencia del nudo trivial, y tiene como representantes los nudos tajada. La existencia de inversos en $\mathcal{C}$ son consecuencias de la siguiente proposici\'on. 

\proposicion{}{ 
Para todo nudo $K\subset S^3$, el nudo $K\# -K$ es la frontera de un disco en $\R^4$.}

\begin{figure}[h]
\centering
\includegraphics[width=0.25\textwidth]{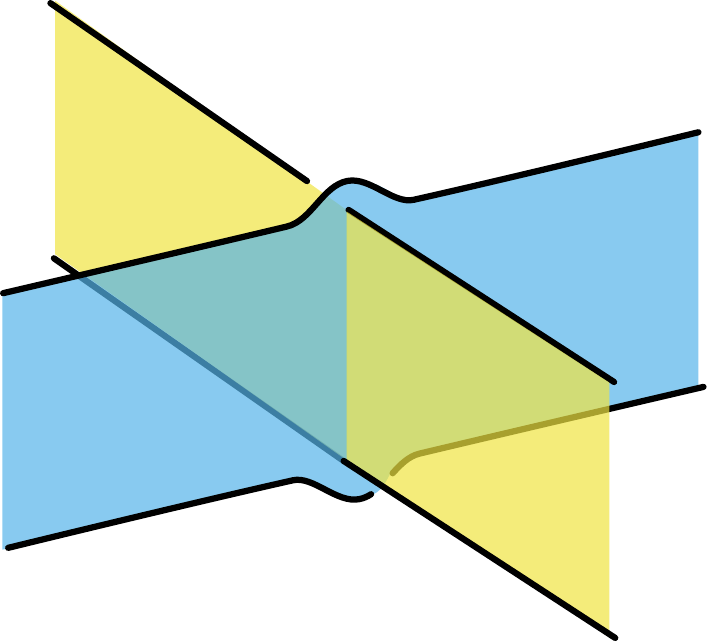}
\caption{Singularidad de list\'on.}\label{ribbon-sing}
\end{figure}

\demostracion{
Coloque una copia de $K$ justo encima del plano $xy$, y una copia de su espejo $m(K)$ justo debajo del plano. N\'otese que lejos de los cruces, los nudos $K$ y $m(K)$ son idénticos, y consisten en la uni\'on disjunta de arcos $\alpha$ ubicados en $\R^3$ en los planos $z=1$ y $z=-1$, digamos. Podemos entonces construir una superficie como el producto de estos arcos con el intervalo $[-1,1]$, y las dos copias de $\alpha$ son precisamente la frontera `horizontal' de \'esta superficie. En contraste, en los cruces del nudo cada hebra de $K$ se une con la hebra correspondiente de $m(K)$ a trav\'es de un `rect\'angulo'. La superficie as\'i construida es la imagen de una inmersi\'on $\psi:[0,1]\times S^1\hookrightarrow\R^3$ cuyo conjunto singular $S=\{x\in \text{Im}(\psi) \mid \psi \text{ no es } 1-1 \}$ tiene una estructura especial: dado que en cada cruce la hebra superior de $K$ corresponde a la hebra inferior de $m(K)$ (y viceversa), cada componente conexa de $S$ es la imagen bajo $\psi$ de un par de intervalos cerrados, uno con puntos l\'imites en la frontera de $[0,1]\times S^1$ y otro completamente contenido en su interior (véase \Cref{ribbon-sing}). Estas singularidades se conocen como singularidades de cinta o list\'on, y se pueden resolver en $\R^4_+$ de la siguiente manera. Sea $I\subset [0,1]\times S^1$ un intervalo completamente contenido en el interior del cilindro $[0,1]\times S^1$ tal que $\psi(I)$ es una componente conexa de $S$. Sea $U$ una vecindad de $I$ en el cilindro. El resultado de incluir el interior de $\psi(U)$ al interior de $\R^4_+$ mientras se deja $\psi(\partial U)$ fija en $\R^3\times\{0\}$ es un cilindro $[0,1]\times S^1$ propiamente encajado en $\R^4_+$. El complemento en el cilindro de un entorno tubular de un arco vertical es el disco deseado. Su frontera es la suma conexa de $K$ y $m(K^r)$, su imagen especular con orientaci\'on opuesta.\footnote{Excelentes figuras que ilustran \'esta construcci\'on se pueden ver en el siguiente enlace \url{https://sketchesoftopology.wordpress.com/2009/04/20/mirrors-and-ribbons/}}.
}

\teorema{Levine \cite{levine1,levine2}}{ Existe un homomorfismo de grupos sobreyectivo $f:\mathcal{C}\to \Z^\infty\oplus (\Z/2)^\infty\oplus (\Z/4)^\infty $.}

De hecho, los elementos de $\mathcal{C}$ cuya imagen bajo $f$ es $ \Z^\infty\oplus (\Z/2)^\infty$ fueron introducidos en \Cref{familias}. Primero, la colecci\'on de nudos $\{T_{2,2k+1}\}_{k=1}^\infty$ de \Cref{toricos} genera un subgrupo de $\mathcal{C}$ isomorfo a $\Z^\infty$, es decir, ninguna combinaci\'on lineal de sus elementos es un nudo `tajada'. C\'alculos similares a los llevados a cabo al final de \Cref{invariantes} dan $\sigma(T_{2,2k+1})=-2k$, y por lo tanto cada nudo $T_{2,2k+1}$ genera un subgrupo isomorfo a $\Z$. Para obtener el subgrupo $\Z^\infty$ se requiere una generalizaci\'on de la signatura de nudos a una funci\'on de signaturas $\sigma_K:S^1\to\Z$.\\

Como consecuencia de la \Cref{anfiquiral} tenemos que existe una isotop\'ia entre el nudo $K(2m,2m)$ (ver la \Cref{double_twist}) y su imagen especular $K(-2m,-2m)$ (con orientaci\'on opuesta). De \'esto se desprende que la suma $K(2m,2m)\#K(2m,2m)$ es un nudo `tajada', o en otras palabras, el nudo $K(2m,2m)$ tiene orden dos en $\mathcal{C}$ para $m\neq 0$. Un an\'alisis cuidadoso de los factores irreducibles del polinomio de Alexander de los elementos de la colecci\'on $\{K(2m,2m)\}_{m=1}^\infty$ permite concluir que \'esta genera un subgrupo isomorfo a $(\Z/2)^\infty$. Los detalles se pueden encontrar en \cite{livingston-naik}.


\bibliographystyle{plain}
\bibliography{thebibliography}

\end{document}